\documentclass{svproc}

\newcommand{\cahiernumber}{02}
\makeatletter
\def\cahierline{\gdef\@cahierline}
\cahierline{\textcolor{gray}{Cahier du GERAD G-\number\year-\cahiernumber}}

\def\ps@mytitlepage{\let\@mkboth\@gobbletwo
   \def\@evenhead{\normalfont\small\rlap{\@cahierline}\hspace{\headlineindent}%
                  \hfil}
   \def\@oddhead{\normalfont\small\hfil\hspace{\headlineindent}%
                 \llap{\@cahierline}}
   \def\chaptermark##1{}%
   \def\sectionmark##1{}%
   \def\subsectionmark##1{}}

\def\ps@myheadings{\let\@mkboth\@gobbletwo
   \def\@oddfoot{\normalfont\scriptsize\ttfamily\@cahierline}
   \def\@evenfoot{\normalfont\scriptsize\ttfamily\@cahierline}
   \def\@evenhead{\normalfont\small\rlap{\thepage}\hspace{\headlineindent}%
    \leftmark\hfil}
    \def\@oddhead{\normalfont\small\scshape\hfil\rightmark\hspace{\headlineindent}%
                 \llap{\thepage}}
   \def\chaptermark##1{}%
   \def\sectionmark##1{}%
   \def\subsectionmark##1{}}
\makeatother

\newlength{\springerbaselineskip}
\usepackage{url}

\usepackage{amsmath}
\usepackage{microtype}
\usepackage{dsfont}
\usepackage{xspace}
\usepackage{xcolor}

\usepackage{longtable}

\usepackage{etoolbox}
\AtBeginEnvironment{longtable}{\scriptsize}

\makeatletter
\def\LT@makecaption#1#2#3{%
  \LT@mcol\LT@cols c{\hbox to\z@{\hss\parbox[t]\LTcapwidth{%
    \sbox\@tempboxa{{\normalsize{\bfseries #1{#2.}} #3}}%
    \ifdim\wd\@tempboxa>\hsize
      \baselineskip \springerbaselineskip
      {\normalsize{\bfseries #1{#2.}} #3}%
    \else
      \hbox to\hsize{\hfil\box\@tempboxa\hfil}%
    \fi
    \endgraf\vskip\baselineskip}%
  \hss}}}
\makeatother

\usepackage[utf8]{inputenc}
\usepackage[documenter]{jlcode}
\usepackage{todonotes}

\newcommand{\R}{\mathds{R}}
\newcommand*{\doi}[1]{DOI: \url{http://dx.doi.org/{#1}}}

\newcommand{\bc}{\begin{center}}
\newcommand{\ec}{\end  {center}}
\newcommand{\bd}{\begin{description}}
\newcommand{\ed}{\end  {description}}
\newcommand{\be}{\begin{enumerate}}
\newcommand{\ee}{\end  {enumerate}}
\newcommand{\bi}{\begin{itemize}}
\newcommand{\ei}{\end  {itemize}}
\newcommand{\bpm}{\begin{pmatrix}}
\newcommand{\epm}{\end  {pmatrix}}
\newcommand{\bs}{\begin{small}}
\newcommand{\es}{\end  {small}}
\newcommand{\bt}{\begin{tabular}}
\newcommand{\et}{\end  {tabular}}
\newcommand{\bft}{\begin{footnotesize}}
\newcommand{\eft}{\end  {footnotesize}}

\newcommand{\pmat}[1]{\begin{pmatrix}#1\end{pmatrix}}

\newcommand{\half}  {{\textstyle\frac12}}
\newcommand{\kp}[1]{_{k+#1}}

\newcommand{\maxim}{\mathop{\mathrm{maximize}}}
\newcommand{\maximize}[1]{{\displaystyle\maxim_{#1}}}
\newcommand{\minim}{\mathop{\mathrm{minimize}}}
\newcommand{\minimize}[1]{{\displaystyle\minim_{#1}}}
\newcommand{\norm}[1]{\|#1\|}
\newcommand{\T}{^T\!}
\newcommand{\z}{\phantom{0}}

\newcommand{\AMPL}{{\small AMPL}\xspace}
\newcommand{\CUTEst}{{\small CUTEst}\xspace}
\newcommand{\GAMS}{{\small GAMS}\xspace}
\newcommand{\JuMP}{{\small JuMP}\xspace}
\newcommand{\LANCELOT}{{\small LANCELOT}\xspace}
\newcommand{\MINOS}{{\small MINOS}\xspace}
\newcommand{\SNOPT}{{\small SNOPT}\xspace}
\newcommand{\IPOPT}{{\small IPOPT}\xspace}
\newcommand{\KNITRO}{{\small KNITRO}\xspace}
\newcommand{\NCL}{{\small NCL}\xspace}
\newcommand{\BCk}{\mbox{BC$_k$}\xspace}
\newcommand{\LCk}{\mbox{LC$_k$}\xspace}
\newcommand{\NCk}{\mbox{NC$_k$}\xspace}
\newcommand{\NCa}{\mbox{NC$_0$}\xspace}
\newcommand{\NCkp}{\mbox{NC$_{k+1}$}\xspace}
\newcommand{\nclnls}{\texttt{ncl\_nls}\xspace}
\newcommand{\knitronls}{\texttt{knitro\_nls}\xspace}

\newcommand{\Problem}[4]{\fbox
   {\begin{tabular*}{0.84\textwidth}
    {@{}l@{\extracolsep{\fill}}l@{\extracolsep{6pt}}%
        l@{\extracolsep{\fill}}c@{}}
      #1 & $\minimize{#2}$ & $#3$ & $ $ \\[5pt]
         & subject to  & $#4$ & $ $
    \end{tabular*}}}

\newcommand{\maxproblem}[4]{\fbox
   {\begin{tabular*}{0.84\textwidth}
    {@{}l@{\extracolsep{\fill}}l@{\extracolsep{6pt}}%
        l@{\extracolsep{\fill}}c@{}}
      #1 & $\maximize{#2}$ & $#3$ & $ $ \\[5pt]
         & subject to  & $#4$ & $ $
    \end{tabular*}}}

\newcommand{\starsymbol}{\ast}
\newcommand{\superstar}{^\starsymbol}

\newcommand{\xstar}{x\superstar}
\newcommand{\ystar}{y\superstar}
\newcommand{\zstar}{z\superstar}
\newcommand{\etastar}{\eta\superstar}
\newcommand{\omegastar}{\omega\superstar}
\newcommand{\rhostar}{\rho\superstar}

\newcommand{\xstark}{x^*_k}
\newcommand{\xstarkm}{x^*_{k-1}}
\newcommand{\rstark}{r^*_k}

\newcommand{\ystark}{y^*_k}
\newcommand{\zstark}{z^*_k}

\begin{document}
\mainmatter              
\pagestyle{myheadings}
\thispagestyle{mytitlepage}

\title{A Julia implementation of Algorithm NCL
    \\ for constrained optimization
    \\ {\normalsize \rm \today}
    \\[-25pt] \null}
\titlerunning{A Julia Implementation of Algorithm NCL}

\author{
  Ding Ma\inst{1} \and Dominique Orban\inst{2} \and Michael A. Saunders\inst{3}%
}
\authorrunning{D. Ma, D. Orban, and M. A. Saunders}
\tocauthor{Ding Ma, Dominique Orban, and Michael A. Saunders}

\institute{Department of Management Science and Department of Marketing,
\\ College of Business, City University of Hong Kong, Hong Kong,
\\ \email{dingma@cityu.edu.hk}, \url{https://www.cb.cityu.edu.hk/staff/dingma}
\and
   GERAD and Department of Mathematics and Industrial Engineering,
\\ Polytechnique Montr\'eal, QC, Canada,
\\ \email{dominique.orban@gerad.ca}, \url{https://dpo.github.io}
\and
   Systems Optimization Laboratory,
   Department of Management Science and Engineering,
   Stanford University,
   Stanford, CA, USA\@.
\\ \email{saunders@stanford.edu}, \url{http://stanford.edu/~saunders}%
}

\maketitle

\begin{abstract}
  Algorithm \NCL is designed for general smooth optimization problems
  where first and second derivatives are available,
  including problems whose constraints may not be linearly independent at
  a solution (i.e., do not satisfy the LICQ).
  It is equivalent to the LANCELOT augmented Lagrangian method,
  reformulated as a short sequence of nonlinearly constrained
  subproblems that can be solved efficiently by \IPOPT and \KNITRO, with
  warm starts on each subproblem.  We give numerical results
  from a Julia implementation of Algorithm \NCL on tax policy models that do not satisfy the LICQ, and on nonlinear least-squares problems and general problems from the \CUTEst test set.
\end{abstract}



\section{Introduction}

Algorithm NCL (nonlinearly constrained augmented Lagrangian \cite{NAOIV18}) is designed for smooth, constrained
optimization problems for which first and second derivatives are
available.  Without loss of generality, we take the problem to be
\[
  \Problem{NCO}{x \in \R^n}{\phi(x)}{c(x) = 0, \quad \ell \le x \le u,}
\]
where $\phi(x)$ is a scalar objective function and $c(x) \in \R^m$ is a vector
of linear or nonlinear constraints.  Inequality constraints are
accommodated by including slack variables within $x$.  We take the
primal and dual solutions to be $(\xstar,\ystar,\zstar)$.  We denote
the objective gradient by $g(x) = \nabla \phi(x) \in \R^n$, and the constraint
Jacobian by $J(x) \in \R^{m \times n}$. The objective
and constraint Hessians are $H_i(x) \in \R^{n \times n}$, $i=0,1,\dots,m$.

If $J(\xstar)$ has full row rank $m$, problem NCO
satisfies the linear independence constraint qualification (LICQ) at
$\xstar$.  Most constrained optimization solvers have difficulty if
NCO does \emph{not} satisfy the LICQ.  An exception is \LANCELOT
\cite{CGT91,CGT92,LANCELOT}.  Algorithm NCL inherits this desirable
property by being \emph{equivalent} to the \LANCELOT algorithm.
Assuming first and second derivatives are available,
Algorithm NCL may be viewed as an \emph{efficient implementation}
of the \LANCELOT algorithm.  Previously we have implemented
Algorithm NCL in \AMPL \cite{AMPL,AMPLbook,NAOIV18}
for tax policy problems \cite{JuddSu2011,NAOIV18} that
could not otherwise be solved.\footnote{Available from https://github.com/optimizers/ncl}
Here we describe our implementation in Julia \cite{bezanson2017julia}
and give results on the tax problems and on a set of nonlinear least-squares problems
from the \CUTEst test set \cite{CUTEst}.

\section{LANCELOT and NCL}

For problem NCO, \LANCELOT implements
what we call a \emph{BCL algorithm} (bound-constrained augmented
Lagrangian algorithm), which solves a sequence of about 10
bound-constrained subproblems
\[
   \Problem{\BCk}{x}
   {L(x,y_k,\rho_k) = \phi(x) - y_k\T c(x) + \half \rho_k \norm{c(x)}^2}
   {\ell \le x \le u}
\]
for $k=0,1,2,\dots$, where $y_k$ is an estimate of the dual variable
associated with $c(x)=0$, and $\rho_k > 0$ is a penalty parameter.
Each subproblem \BCk is solved (approximately) with a decreasing
optimality tolerance $\omega_k$, giving an iterate
$(\xstark,\zstark)$.  If $\norm{c(\xstark)}$ is no larger than a
decreasing feasibility tolerance $\eta_k$, the dual variable is
updated to $y\kp1 = y_k - \rho_k c(\xstark)$.  Otherwise, the penalty
parameter is increased to $\rho\kp1 > \rho_k$.

Optimality is declared if $c(\xstark) \le \eta_k$ and $\eta_k$, $\omega_k$
have already been decreased to specified minimum values $\etastar$,
$\omegastar$.  Infeasibility is declared if $c(\xstark) > \eta_k$ and
$\rho_k$ has already been increased to a specified maximum value $\rhostar$.

If $n$ is large and not many bounds are active at $\xstar$, the \BCk
subproblems have \emph{many degrees of freedom}, and \LANCELOT must
optimize in high-dimension subspaces.  The subproblems are therefore
computationally expensive.  The algorithm in \MINOS \cite{MurS82}
(we call it an LCL algorithm) reduces this expense by including linearizations of the constraints
within its subproblems:
\[
   \Problem{\LCk}{x}
   {L(x,y_k,\rho_k) = \phi(x) - y_k\T c(x) + \half \rho_k \norm{c(x)}^2}
   {c(\xstarkm) + J(\xstarkm)(x - \xstarkm) = 0, \quad \ell \le x \le u.}
\]
The SQP algorithm in \SNOPT \cite{GilMS05} solves subproblems with
the same linearized constraints and a quadratic approximation to the
\LCk objective.  Complications arise for both \MINOS and \SNOPT
if the linearized constraints are infeasible.

Algorithm NCL proceeds in the opposite way by introducing additional
variable $r \equiv - c(x)$ into subproblems \LCk to obtain the NCL
subproblems
\[
   \Problem{\NCk}{x,\,r}
   {\phi(x) + y_k\T r  + \half \rho_k \norm{r}^2}
   {c(x) + r = 0, \quad \ell \le x \le u.}
\]
These subproblems have \emph{nonlinear constraints} and far more degrees of
freedom than the original NCO!  Indeed, the extra variables $r$ make
the subproblems more difficult if they are solved by \MINOS and \SNOPT.
However, the subproblems satisfy the LICQ because of $r$.
Also, interior solvers such as \IPOPT \cite{IPOPT} and
\KNITRO \cite{KNITRO} find $r$ helpful because at each
\emph{interior iteration $p$}
they update the current primal-dual point $(x_p,r_p,\lambda_p)$ by
computing a search direction $(\Delta x,\Delta r,\Delta\lambda)$ from
a linear system of the form
\begin{equation} \label{eq:search}
  \pmat{(H_p+D_p) &                   & J_p^T
     \\          & \!\!\!\!\rho_k I & I
     \\  J_p     &                 I &      }
     \pmat{\Delta x \\ \Delta r \\ \Delta\lambda}
  = -\pmat{  g(x_p) + J_p^T \lambda_p - z_p + w_p \\ y_k + \rho_k r_p + \lambda_p \\ c(x_p) + r_p},
\end{equation}
where $D_p$, \(z_p\) and \(w_p\) are an ill-conditioned positive-definite diagonal matrix and two vectors arising from the interior method, and each Lagrangian Hessian
$H_p = H_0(x_p) - \sum_i (y_k)_i H_i(x_p)$ may be altered to be more
positive definite.  Direct methods for solving each sparse system
\eqref{eq:search} are affected very little by the higher dimension
caused by $r$, and they benefit significantly from $\pmat{J_p & I}$
always having full row rank.

If an optimal solution for \NCk is $(\xstark,\rstark,\ystark,\zstark)$
and the feasibility and optimality tolerances have decreased to their
minimum values $\eta^*$ and $\omega^*$, a natural stopping condition for
Algorithm NCL is $\norm{\rstark}_\infty \le \eta^*$, because the major
iterations drive $r$ toward zero and we see that if $r=0$, subproblem \NCk
is equivalent to the original problem NCO.

We have found that Algorithm NCL is successful in practice because
\begin{itemize}
\item there are only about 10 major iterations ($k = 1, 2, \dots, 10$);
\item the search-direction computation \eqref{eq:search} for interior
  solvers is more stable than if the solvers are applied to NCO
  directly;
\item \IPOPT and \KNITRO have run-time options that facilitate
  \emph{warm starts} for each subproblem \NCk, $k>1$.
\end{itemize}

\section{Optimal tax policy problems}

The above observations were confirmed by our \AMPL implementation of
Algorithm NCL in solving some large problems modeling taxation policy \cite{JuddSu2011,NAOIV18,NCLweb}.
The problems have very many nonlinear inequality constraints
$c(x) \ge 0$ in relatively few variables.
They have the form
\[
   \maxproblem{TAX}{c,\,y}
      {\hspace*{33pt} \sum_i \lambda_i U^i(c_i,y_i)}
      {\begin{array}[t]{rcl}
          U^i(c_i,y_i) - U^i(c_j,y_j) \ge 0\phantom, &\text{ \ for all } i,j
       \\ \lambda^T(y - c)            \ge 0\phantom, &
       \\          c,\ y              \ge 0,         &
       \end{array}}
\]
where $c_i$ and $y_i$ are the consumption and income of taxpayer $i$,
and $\lambda$ is a vector of positive weights.%
\footnote{In this section, $(c,y,\lambda)$ refer to problem TAX, not the variables
in Algorithm NCL.}
The utility functions $U^i(c_i,y_i)$ are each of the form
\[
   U(c,y) = \frac{(c - \alpha)^{1 - 1/\gamma}}{1 - 1/\gamma}
   - \psi \frac{(y/w)^{1/\eta + 1}}{1/\eta + 1},
\]
where $w$ is the wage rate and $\alpha$, $\gamma$,
$\psi$ and $\eta$ are taxpayer heterogeneities.
More precisely, the utility functions are of the form
\begin{align*}
   U^{i,j,k,g,h}(c_{p,q,r,s,t},y_{p,q,r,s,t}) =
   \frac{(c_{p,q,r,s,t}-\alpha_k)^{1-1/\gamma_h}}{1-1/\gamma_h}
   - \psi_g \frac{(y_{p,q,r,s,t}/w_i)^{1/\eta_j + 1}}{1/\eta_j + 1},
   \label{E:utility5}
\end{align*}
where $(i,j,k,g,h)$ and $(p,q,r,s,t)$ run over
$na$ wage types,
$nb$ elasticities of labor supply,
$nc$ basic need types,
$nd$ levels of distaste for work, and
$ne$ elasticities of demand for consumption,
with $na$, $nb$, $nc$, $nd$, $ne$ determining the size of the problem,
namely $m=T(T-1)$ nonlinear constraints, $n=2T$ variables,
with $T := na \times nb \times nc \times nd \times ne$.

To achieve reliability, we found it necessary to extend the \AMPL model's definition of $U(c,y)$ to be a piecewise-continuous function that accommodates negative values of $(c-\alpha)$.

At a solution, a large proportion of
the constraints are essentially active.
The failure of LICQ causes
numerical difficulties for \MINOS, \SNOPT, and \IPOPT.  \LANCELOT is
more able to find a solution, except it is very slow on each subproblem
\NCk.  For example, on the smallest problem of
Table~\ref{tab:IPOPT-KNITRO} with 32220 constraints and 360 variables,
\LANCELOT running on NEOS \cite{NEOS} timed-out at a near-optimal point on the 11th major iteration after 8 hours of CPU.

Note that when the constraints of NCO are inequalities $c(x) \ge 0$
as in problem TAX,
the constraints of subproblem \NCk become inequalities $c(x) + r \ge 0$
(and similarly for mixtures of equalities and inequalities).  The
inequalities mean ``more free variables'' (more variables that are
not on a bound).  This increases the problem difficulty for \MINOS
and \SNOPT, but has only a positive effect on the interior solvers.

\begin{table}[t]   
\caption{Run-time options for warm-starting \IPOPT and \KNITRO
         on subproblem \NCk.}
\label{tab:options}
\bc
\small
\bt{ll@{\quad}l}
             & IPOPT                               & KNITRO
\\ \hline
   $k=1$     &                                     & {\tt algorithm=1}
\\ $k\ge  2$ & {\tt warm\_start\_init\_point=yes}  & {\tt bar\_directinterval=0}
\\           &                                     & {\tt bar\_initpt=2}
\\           &                                     & {\tt bar\_murule=1}
\\ $k = 2,3$ & {\tt mu\_init=1e-4}                 & {\tt bar\_initmu=1e-4}
\\           &                                     & {\tt bar\_slackboundpush=1e-4}
\\ $k = 4,5$ & {\tt mu\_init=1e-5}                 & {\tt bar\_initmu=1e-5}
\\           &                                     & {\tt bar\_slackboundpush=1e-5}
\\ $k = 6,7$ & {\tt mu\_init=1e-6}                 & {\tt bar\_initmu=1e-6}
\\           &                                     & {\tt bar\_slackboundpush=1e-6}
\\ $k = 8,9$ & {\tt mu\_init=1e-7}                 & {\tt bar\_initmu=1e-7}
\\           &                                     & {\tt bar\_slackboundpush=1e-7}
\\ $k\ge 10$ & {\tt mu\_init=1e-8}                 & {\tt bar\_initmu=1e-8}
\\           &                                     & {\tt bar\_slackboundpush=1e-8}
\\ \hline
\et
\ec
\vspace*{-5pt}
\end{table}

\begin{table}[t]  
  \caption{Solution of tax problems of increasing dimension using IPOPT
    and KNITRO on the original problem (cold starts) and the AMPL
    implementation of Algorithm NCL with IPOPT or KNITRO as subproblem
    solvers (warm-starting with the options in Table~\ref{tab:options}).
    The problem size increases with a problem parameter $na$.  Other
    problem parameters are fixed at $nb=nc=3$, $nd=ne=2$.  There are $m$
    nonlinear inequality constraints and $n$ variables.  For IPOPT, $>$
    indicates optimality was not achieved.}
  \label{tab:IPOPT-KNITRO}
\small

\bc
\bt{ccc|rr|rr|rr|rr}
     &        &      & \multicolumn{2}{c|}{IPOPT} & \multicolumn{2}{c|}{KNITRO}
                     & \multicolumn{2}{c}{{NCL/IPOPT}} & \multicolumn{2}{|c}{NCL/KNITRO}
\\ \hline
 $na$&    $m$ &  $n$ &      itns &       time    &    itns &  time & {itns} & {time}  & itns &  time
\\ \hline
   5 &  32220 &  360 &       449 &       217     &     168 &    53 & { 322} & {  146} &  339 &    63 
\\ 9 & 104652 &  648 &      $>98$&     $>360$    &     928 &   825 & { 655} & { 1023} &  307 &   239 
\\11 & 156420 &  792 &      $>87$&     $>600$    &    2769 &  4117 & { 727} & { 1679} &  383 &   420
\\17 & 373933 & 1224 &           &               &    2598 & 11447 & {1021} & { 6347} &  486 &  1200
\\21 & 570780 & 1512 &           &               &         &       & {1761} & {17218} &\z712 &  2880
\\ \hline
\et%
\ec
\vspace*{-5pt}
\end{table}

In wishing to improve the efficiency of Algorithm NCL on larger tax
problems, we found it possible to warm-start \IPOPT and \KNITRO on
each \NCk subproblem ($k>1$) by setting the run-time options shown
in Table~\ref{tab:options}.
These options were used by \NCL/\IPOPT and \NCL/\KNITRO
to obtain the results in Table~\ref{tab:IPOPT-KNITRO}.
We see that \NCL/\IPOPT performed significantly better than \IPOPT itself,
and similarly for \NCL/\KNITRO compared to \KNITRO.
The feasibility and optimality tolerances $\eta_k$, $\omega_k$
were fixed at $\eta^* = \omega^* = \hbox{\tt 1e-6}$ for all $k$.
Our Julia implementation saves computation by starting
with larger $\eta_k$, $\omega_k$ and reducing them toward $\eta^*$, $\omega^*$
as in \LANCELOT.


\newpage

\section{Julia implementation}

Modeling languages such as \AMPL and \GAMS are domain-specific languages, as opposed to full-fledged, general-purpose programming languages like C or Java.
In the terminology of Bentley \cite{bentley-1986}, they are \emph{little languages}.
As such, they have understandable, yet very real limitations, that make it difficult, impractical, and perhaps even impossible, to implement an algorithm such as Algorithm \NCL in a sufficiently generic manner so that it may be applied to arbitrary problems.
Indeed, our \AMPL implementation of Algorithm \NCL is specific to the optimal tax policy problems, and it would be difficult to generalize it to other problems.
One of the main motivations for implementing Algorithm \NCL in a language such as Julia is to be able to solve a greater variety of optimization problems.

We now describe the key features of our Julia implementation of Algorithm \NCL and show that it solves examples of the same tax problems more efficiently.
We then give results on a set of nonlinear least-squares problems from the \CUTEst test set to indicate that Algorithm \NCL is a reliable solver for such problems where first and second derivatives are available for the interior solvers used at each major iteration.
To date, this means that Algorithm \NCL is
effective for optimization problems modeled in \AMPL, \GAMS, and \CUTEst.  (We have not made an implementation in \GAMS \cite{GAMS}, but it would be possible to build a major-iteration loop around calls to \IPOPT or \KNITRO in the way that we did for \AMPL \cite{NAOIV18}.)

\subsection{Key features}

The main advantage of a Julia implementation over our original \AMPL implementation is that we may take full advantage of our Julia software suite for optimization, hosted under the \emph{JuliaSmoothOptimizers} (JSO) organization \cite{orban-siqueira-jso-2020}.
Our suite provides a general consistent API for solvers to interact with models by providing flexible data types to represent the objective and constraint functions, to evaluate their derivatives, to examine bounds on the variables, to add slack variables transparently, and to provide essentially any information that a solver might request from a model.
Thanks to interfaces to modeling languages such as \AMPL, \CUTEst and \JuMP \cite{jump}, solvers in JSO may be written without regard for the language in which the model was written.

The modules from our suite that are particularly useful in the context of our implementation of Algorithm NCL are the following.
\begin{itemize}
    \item NLPModels \cite{orban-siqueira-nlpmodels-2020} is the main modeling package that defines the API on which solvers can rely to interact with models.
    Models are represented as instances of a data type deriving from the base type \texttt{AbstractNLPModel}, and solvers can evaluate the objective value by calling the \texttt{obj()} method, the gradient vector by calling the \texttt{grad()} method, and so forth.
    The main advantage of the consistent API provided by NLPModels is that solvers need not worry about the provenance of models.
    Other modules ensure communication between modeling languages such as \AMPL, \CUTEst or \JuMP, and NLPModels.
    \item AmplNLReader \cite{orban-siqueira-amplnlreader-2020} is one such module, and, as the name indicates, allows a solver written in Julia to interact with a model written in \AMPL.
    The communication is made possible by the \AMPL Solver Library (ASL)\footnote{\url{http://www.netlib.org/ampl/solvers}}, which requires that the model be decoded as an {\tt nl} file.
    \item NLPModelsIpopt \cite{orban-siqueira-nlpmodelsipopt-2020} is a thin translation layer between the low-level Julia interface to \IPOPT provided by the IPOPT.jl package\footnote{\url{https://github.com/jump-dev/Ipopt.jl}} and NLPModels, and lets users solve any problem conforming to the NLPModels API with \IPOPT.
    \item NLPModelsKnitro \cite{orban-siqueira-nlpmodelsknitro-2020} is similar to NLPModelsIpopt, but lets users solve problems with \KNITRO via the low-level interface provided by KNITRO.jl\footnote{\url{https://github.com/jump-dev/KNITRO.jl}}.
\end{itemize}

Julia is a convenient language built on top of state-of-the-art infrastructure underlying modern compilers such as Clang.
Julia may be used as an interactive language for exploratory work in a read-eval-print loop similar to Matlab.
However, Julia functions are transparently translated to low-level code and compiled the first time they are called.
The net result is efficient compiled code whose efficiency rivals that of binaries generated from standard compiled languages such as C and Fortran.
Though this last feature is not particularly important in the context of Algorithm NCL because the compiled solvers \IPOPT and \KNITRO perform all the work, it is paramount when implementing pure Julia optimization solvers.

\subsection{Implementation and solver features}

The Julia implementation of Algorithm NCL, named NCL.jl \cite{orban-personnaz-ncl-2020}, is in two parts.
The first part defines a data type \texttt{NCLModel} that derives from the basic data type \texttt{AbstractNLPModel} mentioned earlier and represents subproblem \NCk.
An \texttt{NCLModel} is a wrapper around the underlying problem NCO in which the current values of \(\rho_k\) and \(r\) can be updated efficiently.
The second part is the solver itself, each iteration of which consists of a call to \IPOPT or \KNITRO, and parameter updates.
The solver takes an \texttt{NCLModel} as input.
If the input problem is not an \texttt{NCLModel}, it is first converted into one.
Parameters are initialized as
\[
    \eta_0 = 10, \quad
    \omega_0 = 10, \quad
    \rho_0 = 100, \quad
    \mu_0 = 0.1,
\]
where \(\mu_0\) is the initial barrier parameter for \IPOPT or \KNITRO.
The initial values of \(x\) are those defined in the underlying model if any, or zero otherwise.
We initialize \(r\) to zero and \(y\) to the vector of ones.
When the subproblem solver returns with \NCk solution $(\xstark,\rstark,\ystark,\zstark)$, we check whether \(\norm{\rstark} \leq \max(\eta_k, \, \eta_*)\).
If so, we decide that good progress has been made toward feasibility and update
\[
    y_{k+1} = y_k - \rho_k \rstark, \quad
    \eta_{k+1} = \eta_k / 10, \quad
    \omega_{k+1} = \omega_k / 10, \quad
    \rho_{k+1} = \rho_k,
\]
where this definition of $y_{k+1}$ is the first-order update of the multipliers.
Otherwise, we keep most things the same but increase the penalty parameter:
\[
    y_{k+1} = y_k, \quad
    \eta_{k+1} = \eta_k, \quad
    \omega_{k+1} = \omega_k, \quad
    \rho_{k+1} = \min(10 \rho_k, \, \rho^*),
\]
where \(\rhostar > 0\) is the threshold beyond which the user is alerted that the problem may be infeasible.
In our implementation, we use \(\rhostar = 10^{12}\).

Note that updating the multipliers based on \(\norm{\rstark}\) instead of \(\norm{c(\xstark)}\) is a departure from the classical augmented-Lagrangian update.
From the optimality conditions for \NCk we can prove that the first-order update is equivalent to choosing $y_{k+1} = \ystark$ when \NCk is solved accurately.
We still have a choice between the two updates because we use low accuracy for the early \NCk.  We could also ``trim'' $y_{k+1}$ (i.e., for inequality constraints $c_i(x) + r_i \ge 0$ or $\le 0$, set components of $y_{k+1}$ with non-optimal sign to zero).  These are topics for future research.

With \IPOPT as subproblem solver, we warm-start subproblem \NCkp with
the options in Table~\ref{tab:options} and
$(\ystark, \zstark)$ as initial values for the Lagrange multipliers.
With \KNITRO as subproblem solver, $(\ystark, \zstark)$ as starting point did not help or harm \KNITRO significantly.
We allowed \KNITRO to determine its own initial multipliers, and it proved to be
significantly more reliable than \IPOPT in solving the \NCk subproblems for the optimal tax policy problems.
In the next sections, Algorithm NCL means our Julia implementation with \KNITRO as subproblem solver.

\subsection{Results with Julia/NCL on the tax policy problems}

\AMPL models of the optimal tax policy problems were input to
the Julia implementation of Algorithm NCL.
The notation 1D, 2D, 3D, 4D, 5D refers to problem parameters
$na$, $nb$, $nc$, $nd$, $ne$ that define the utility function
appearing in the objective and constraints.
The subproblem solver was \KNITRO 12 \cite{KNITRO}.

Tables \ref{tab:tax1D}--\ref{tab:tax5D} illustrate that, as with our
\AMPL implementation of Algorithm NCL, about
10 major iterations are needed independent of the problem size.
(The problems have increasing numbers of variables and greatly increasing
numbers of nonlinear inequality constraints.)
In each iteration log,

\smallskip

\begin{quote}
{\tt outer} and {\tt inner} refer to the
NCL major iteration number $k$ and the total number of \KNITRO iterations
for subproblems \NCk;

\smallskip

{\tt NCL obj} is the augmented Lagrangian objective value, which
converges to the objective value for the model;

\smallskip

$\eta$ and $\omega$ show the \KNITRO feasibility and optimality tolerances $\eta_k$ and $\omega_k$ decreasing from $10^{-2}$ to $10^{-6}$;

\smallskip

$\norm{\nabla L}$ is the size of the augmented Lagrangian gradient,
namely $\norm{g(\xstark) - J(\xstark)^T y_{k+1}}$ (a measure of the dual infeasibility at the end of major iteration $k$);

\smallskip

$\rho$ is the penalty parameter $\rho_k$;

\smallskip

$\mu$ {\tt init} is the initial value of \KNITRO's barrier parameter;

\smallskip

$\norm{x}$ is the size of the primal variable $\xstark$ at the
(approximate) solution of \NCk;

\smallskip

$\norm{y}$ is the size of the corresponding dual variable $\ystark$;

\smallskip

{\tt time} is the number of seconds to solve \NCk.
\end{quote}

\smallskip
\noindent
We see from the decreasing {\tt inner} iteration counts that \KNITRO was able to warm-start each subproblem, and from the decreasing $\norm{r}$ and $\norm{\nabla L}$ values that it is sufficient to solve the early subproblems with low (but steadily increasing) accuracy.

\begin{table}[htbp]  
  \caption{Tax1D problem with realistic data.
  NCL with KNITRO solving subproblems.}
  \label{tab:tax1D}
  \centering
  \vspace*{-20pt}
  \begin{lstlisting}
julia> using NCL

julia> using AmplNLReader

julia> tax1D = AmplModel("data/tax1D")
Maximization problem data/tax1D
nvar = 24, ncon = 133 (1 linear)

julia> NCLSolve(tax1D, outlev=0)

outer inner   NCL obj     ‖r‖       η    ‖∇L‖       ω       ρ  μ init     ‖y‖     ‖x‖ time
    1     5 -8.00e+02 9.7e-02 1.0e-02 7.6e-03 1.0e-02 1.0e+02 1.0e-01 1.0e+00 2.0e+02 0.13
    2    12 -7.89e+02 4.2e-02 1.0e-02 4.3e-03 1.0e-02 1.0e+03 1.0e-03 1.0e+00 1.9e+02 0.00
    3     7 -7.83e+02 5.7e-03 1.0e-02 1.0e-03 1.0e-02 1.0e+04 1.0e-03 1.0e+00 1.9e+02 0.00
    4     3 -7.82e+02 1.3e-04 1.0e-03 1.0e-05 1.0e-03 1.0e+04 1.0e-05 5.8e+01 1.9e+02 0.00
    5     2 -7.82e+02 2.3e-06 1.0e-04 1.0e-05 1.0e-04 1.0e+04 1.0e-05 5.9e+01 1.9e+02 0.00
    6     2 -7.82e+02 9.3e-08 1.0e-05 1.0e-06 1.0e-05 1.0e+04 1.0e-06 5.9e+01 1.9e+02 0.00
    7     2 -7.82e+02 7.7e-09 1.0e-06 1.0e-08 1.0e-06 1.0e+04 1.0e-06 5.9e+01 1.9e+02 0.00
  \end{lstlisting}
\end{table}

\begin{table}[htbp]  
  \caption{Tax2D problem.  NCL with KNITRO solving subproblems.}
  \label{tab:tax2D}
  \centering
  \vspace*{-20pt}
  \begin{lstlisting}
julia> tax2D = AmplModel("data/tax2D")
Maximization problem data/tax2D
nvar = 120, ncon = 3541 (1 linear)

julia> NCLSolve(tax2D, outlev=0)
outer inner   NCL obj     ‖r‖       η    ‖∇L‖       ω       ρ  μ init     ‖y‖     ‖x‖ time
    1    16 -4.35e+03 6.1e-02 1.0e-02 4.2e-03 1.0e-02 1.0e+02 1.0e-01 1.0e+00 4.0e+02 0.15
    2    15 -4.31e+03 2.5e-02 1.0e-02 2.7e-04 1.0e-02 1.0e+03 1.0e-03 1.0e+00 4.0e+02 0.13
    3    16 -4.29e+03 7.8e-03 1.0e-02 3.5e-04 1.0e-02 1.0e+04 1.0e-03 1.0e+00 4.0e+02 0.16
    4    15 -4.28e+03 5.1e-03 1.0e-03 1.0e-05 1.0e-03 1.0e+04 1.0e-05 7.9e+01 4.0e+02 0.14
    5    32 -4.28e+03 1.2e-03 1.0e-03 1.0e-05 1.0e-03 1.0e+05 1.0e-05 7.9e+01 4.0e+02 0.32
    6    12 -4.28e+03 1.5e-04 1.0e-03 1.5e-05 1.0e-03 1.0e+06 1.0e-06 7.9e+01 4.0e+02 0.15
    7     4 -4.28e+03 1.8e-05 1.0e-04 2.7e-06 1.0e-04 1.0e+06 1.0e-06 2.0e+02 4.0e+02 0.06
    8     4 -4.28e+03 1.2e-06 1.0e-05 1.3e-07 1.0e-05 1.0e+06 1.0e-07 2.0e+02 4.0e+02 0.05
    9     3 -4.28e+03 3.5e-07 1.0e-06 1.0e-07 1.0e-06 1.0e+06 1.0e-07 2.0e+02 4.0e+02 0.05
   \end{lstlisting}
\end{table}

\begin{table}[htbp]  
  \caption{Tax3D problem.  NCL with KNITRO solving subproblems.}
  \label{tab:tax3D}
  \centering
  \vspace*{-20pt}
  \begin{lstlisting}
julia> pTax3D = AmplModel("data/pTax3D")
Maximization problem data/pTax3D
nvar = 216, ncon = 11557 (1 linear)

julia> NCLSolve(pTax3D, outlev=0)
outer inner   NCL obj     ‖r‖       η    ‖∇L‖       ω       ρ  μ init     ‖y‖     ‖x‖ time
    1     9 -6.97e+03 4.5e-02 1.0e-02 9.1e-03 1.0e-02 1.0e+02 1.0e-01 1.0e+00 5.7e+02 0.54
    2    18 -6.87e+03 1.7e-02 1.0e-02 2.4e-04 1.0e-02 1.0e+03 1.0e-03 1.0e+00 5.6e+02 0.99
    3    16 -6.83e+03 7.8e-03 1.0e-02 1.7e-03 1.0e-02 1.0e+04 1.0e-03 1.0e+00 5.7e+02 1.01
    4    17 -6.81e+03 5.2e-03 1.0e-03 1.5e-05 1.0e-03 1.0e+04 1.0e-05 7.9e+01 5.6e+02 0.99
    5    54 -6.80e+03 2.6e-03 1.0e-03 1.2e-05 1.0e-03 1.0e+05 1.0e-05 7.9e+01 5.6e+02 3.15
    6    22 -6.80e+03 4.5e-04 1.0e-03 8.0e-05 1.0e-03 1.0e+06 1.0e-06 7.9e+01 5.6e+02 1.30
    7     9 -6.80e+03 1.1e-04 1.0e-04 1.0e-06 1.0e-04 1.0e+06 1.0e-06 5.2e+02 5.6e+02 0.56
    8     8 -6.80e+03 1.1e-05 1.0e-04 1.1e-07 1.0e-04 1.0e+07 1.0e-07 5.2e+02 5.6e+02 0.49
    9     5 -6.80e+03 1.1e-06 1.0e-05 1.0e-07 1.0e-05 1.0e+07 1.0e-07 5.3e+02 5.6e+02 0.32
   10     3 -6.80e+03 8.9e-08 1.0e-06 1.0e-08 1.0e-06 1.0e+07 1.0e-08 5.3e+02 5.6e+02 0.22
   \end{lstlisting}
\end{table}

\begin{table}[t]  
  \caption{Tax4D problem.  NCL with KNITRO solving subproblems.}
  \label{tab:tax4D}
  \centering
  \vspace*{-20pt}
  \begin{lstlisting}
julia> pTax4D = AmplModel("data/pTax4D")
Minimization problem data/pTax4D
nvar = 432, ncon = 46441 (1 linear)

julia> NCLSolve(pTax4D, outlev=0)
outer inner   NCL obj     ‖r‖       η    ‖∇L‖       ω       ρ  μ init     ‖y‖     ‖x‖  time
    1    12 -1.34e+04 3.3e-02 1.0e-02 5.4e-03 1.0e-02 1.0e+02 1.0e-01 1.0e+00 7.2e+02  3.38
    2    12 -1.31e+04 1.3e-02 1.0e-02 4.0e-03 1.0e-02 1.0e+03 1.0e-03 1.0e+00 7.2e+02  3.23
    3    15 -1.30e+04 5.1e-03 1.0e-02 2.1e-04 1.0e-02 1.0e+04 1.0e-03 1.0e+00 7.1e+02  3.86
    4    31 -1.30e+04 3.2e-03 1.0e-03 1.3e-05 1.0e-03 1.0e+04 1.0e-05 5.2e+01 7.0e+02  7.95
    5    37 -1.30e+04 1.8e-03 1.0e-03 1.2e-05 1.0e-03 1.0e+05 1.0e-05 5.2e+01 7.0e+02  9.89
    6    44 -1.29e+04 5.0e-04 1.0e-03 1.1e-06 1.0e-03 1.0e+06 1.0e-06 5.2e+01 7.0e+02 11.93
    7    16 -1.29e+04 2.6e-04 1.0e-04 1.2e-05 1.0e-04 1.0e+06 1.0e-06 5.3e+02 7.0e+02  3.74
    8    30 -1.29e+04 4.4e-05 1.0e-04 1.2e-07 1.0e-04 1.0e+07 1.0e-07 5.3e+02 7.0e+02  8.15
    9     9 -1.29e+04 2.3e-05 1.0e-05 1.2e-07 1.0e-05 1.0e+07 1.0e-07 8.2e+02 7.0e+02  2.49
   10    11 -1.29e+04 3.8e-06 1.0e-05 1.0e-08 1.0e-05 1.0e+08 1.0e-08 8.2e+02 7.0e+02  3.09
   11     6 -1.29e+04 1.7e-07 1.0e-06 1.3e-08 1.0e-06 1.0e+08 1.0e-08 9.4e+02 7.0e+02  1.74
   \end{lstlisting}
\end{table}

\begin{table}[t]  
  \caption{Tax5D problem.  NCL with KNITRO solving subproblems.}
  \label{tab:tax5D}
  \centering
  \vspace*{-20pt}
  \begin{lstlisting}
julia> pTax5D = AmplModel("data/pTax5D")
Minimization problem data/pTax5D
nvar = 864, ncon = 186193 (1 linear)

julia> NCLSolve(pTax5D, outlev=0)
outer inner   NCL obj     ‖r‖       η    ‖∇L‖       ω       ρ  μ init     ‖y‖     ‖x‖   time
    1    64 -1.76e+05 2.0e-01 1.0e-02 2.3e-03 1.0e-02 1.0e+02 1.0e-01 1.0e+00 1.1e+04  80.43
    2    29 -1.74e+05 4.9e-02 1.0e-02 1.2e-03 1.0e-02 1.0e+03 1.0e-03 1.0e+00 1.1e+04  35.02
    3    23 -1.74e+05 1.6e-02 1.0e-02 1.0e-03 1.0e-02 1.0e+04 1.0e-03 1.0e+00 1.1e+04  28.96
    4    46 -1.74e+05 4.1e-03 1.0e-02 3.6e-05 1.0e-02 1.0e+05 1.0e-05 1.0e+00 1.1e+04  54.50
    5    41 -1.74e+05 2.8e-03 1.0e-03 1.7e-05 1.0e-03 1.0e+05 1.0e-05 4.1e+02 1.1e+04  52.72
    6    28 -1.74e+05 6.1e-04 1.0e-03 1.0e-06 1.0e-03 1.0e+06 1.0e-06 4.1e+02 1.1e+04  34.38
    7    13 -1.74e+05 2.1e-04 1.0e-04 1.4e-06 1.0e-04 1.0e+06 1.0e-06 1.0e+03 1.1e+04  14.81
    8    12 -1.74e+05 5.3e-05 1.0e-04 1.2e-07 1.0e-04 1.0e+07 1.0e-07 1.0e+03 1.1e+04  14.80
    9     7 -1.74e+05 4.5e-06 1.0e-05 1.0e-07 1.0e-05 1.0e+07 1.0e-07 1.0e+03 1.1e+04   9.49
   10     5 -1.74e+05 8.0e-07 1.0e-06 1.2e-08 1.0e-06 1.0e+07 1.0e-08 1.0e+03 1.1e+04   7.02
   \end{lstlisting}
\end{table}


\subsection{Results with Julia/NCL on CUTEst test set}

Our Julia module CUTEst.jl \cite{orban-siqueira-cutest-2020} provides an interface with the \CUTEst \cite{CUTEst} environment and problem collection.
Its main feature is to let users instantiate problems from \CUTEst using the \texttt{CUTEstModel} constructor so they can be manipulated transparently or passed to a solver like any other \texttt{NLPModel}.

On a set of \(166\) constrained problems with at least \(100\) variables whose constraints are all nonlinear, \KNITRO solves \(147\) and \NCL solves \(126\).
Although our simple implementation of \NCL is not competitive with plain \KNITRO in general, it does solve a few problems on which \KNITRO fails.
Those are summarized in Tables~\ref{tab:knitro-cutest} and~\ref{tab:ncl-cutest}.
The above results suggest that \NCL's strength might reside in solving difficult problems (rather than being the fastest), and that more research is needed to improve its efficiency.


\begin{table}.  
\caption{%
    \label{tab:knitro-cutest}
    \KNITRO results on CUTEst constrained problems
    (a subset that failed).
}
\scriptsize
\hbox to \textwidth{\hss
\bt{lrrrccrrrrrrrl}
\hline
   name & nvar & ncon & \(f\)\hspace*{14pt} & \(\|\nabla L\|_2\) & \(\|c\|_2\) & \(t\)\hspace*{7pt} & iter & \#\(f\) & \#\(\nabla f\) & \#\(c\) & \#\(\nabla c\) & \#\(\nabla^2 L\) & status
\\ \hline
CATENARY & \(  3003\) & \(  1000\) & \(-2.01\)e\(+10\) & \( 4.9\)e\(+00\) & \( 2.0\)e\(+09\) & \( 18.50\) & \(  2000\) & \(  7835\) & \(  2002\) & \(  7835\) & \(  2002\) & \(  2000\) & max\_iter \\
COSHFUN & \(  6001\) & \(  2000\) & \(-9.82\)e\(+17\) & \( 5.0\)e\(-01\) & \( 0.0\)e\(+00\) & \( 19.40\) & \(  2000\) & \(  2001\) & \(  2001\) & \(  2001\) & \(  2001\) & \(  2000\) & max\_iter \\
DRCAVTY1 & \(  4489\) & \(  3969\) & \( 0.00\)e\(+00\) & \( 0.0\)e\(+00\) & \( 2.2\)e\(-03\) & \( 456.00\) & \(  2000\) & \(  9191\) & \(  2002\) & \(  9191\) & \(  2002\) & \(  2000\) & max\_iter \\
EG3 & \( 10001\) & \( 20000\) & \( 5.11\)e\(+05\) & \( 2.0\)e\(+03\) & \( 3.3\)e\(-01\) & \( 7.18\) & \(    51\) & \(    55\) & \(    52\) & \(    55\) & \(    52\) & \(    52\) & infeasible \\
JUNKTURN & \( 10010\) & \(  7000\) & \( 1.78\)e\(-03\) & \( 1.0\)e\(-02\) & \( 6.4\)e\(-07\) & \( 123.00\) & \(  1913\) & \( 15051\) & \(  1915\) & \( 15051\) & \(  1915\) & \(  1914\) & unknown \\
LUKVLE11 & \(  9998\) & \(  6664\) & \( 5.12\)e\(+04\) & \( 5.1\)e\(+02\) & \( 4.1\)e\(-01\) & \( 86.50\) & \(  2000\) & \(  6945\) & \(  2001\) & \(  6945\) & \(  2001\) & \(  2000\) & max\_iter \\
LUKVLE17 & \(  9997\) & \(  7497\) & \( 3.22\)e\(+04\) & \( 1.8\)e\(-02\) & \( 9.9\)e\(-07\) & \( 47.00\) & \(  2000\) & \(  3190\) & \(  2001\) & \(  3190\) & \(  2001\) & \(  2000\) & max\_iter \\
LUKVLE18 & \(  9997\) & \(  7497\) & \( 1.12\)e\(+04\) & \( 4.0\)e\(+01\) & \( 2.5\)e\(-08\) & \( 83.90\) & \(  2000\) & \(  4190\) & \(  2001\) & \(  4190\) & \(  2001\) & \(  2000\) & max\_iter \\
ORTHRDS2 & \(  5003\) & \(  2500\) & \( 7.62\)e\(+02\) & \( 3.7\)e\(-01\) & \( 5.0\)e\(-13\) & \( 0.70\) & \(    42\) & \(    92\) & \(    43\) & \(    92\) & \(    43\) & \(    43\) & unknown
\\ \hline
\et
\hss}
\end{table}

\begin{table}  
\caption{%
    \label{tab:ncl-cutest}
    \NCL results on the same problems (all successful).
}
\scriptsize
\hbox to \textwidth{\hss
\bt{lrrrccrrrrrrrl}
\hline
name & nvar & ncon & \(f\)\hspace*{14pt} & \(\|\nabla L\|_2\) & \(\|c\|_2\) & \(t\)\hspace*{7pt} & iter & \#\(f\) & \#\(\nabla f\) & \#\(c\) & \#\(\nabla c\) & \#\(\nabla^2 L\) & status
\\\hline
CATENARY & \(  3003\) & \(  1000\) & \(-2.10\)e\(+06\) & \( 1.54\)e\(-09\) & \( 1.00\)e\(-07\) & \( 1.76\) & \(   183\) & \(   430\) & \(   195\) & \(   430\) & \(   206\) & \(   194\) & first\_order \\
COSHFUN & \(  6001\) & \(  2000\) & \(-7.81\)e\(-01\) & \( 9.02\)e\(-07\) & \( 1.00\)e\(-07\) & \( 6.46\) & \(   328\) & \(  1712\) & \(   337\) & \(  1712\) & \(   346\) & \(   337\) & first\_order \\
DRCAVTY1 & \(  4489\) & \(  3969\) & \( 0.00\)e\(+00\) & \( 1.19\)e\(-08\) & \( 1.00\)e\(-06\) & \( 29.30\) & \(   222\) & \(   344\) & \(   233\) & \(   344\) & \(   243\) & \(   232\) & first\_order \\
EG3 & \( 10001\) & \( 20000\) & \( 1.94\)e\(-07\) & \( 1.00\)e\(-08\) & \( 1.00\)e\(-07\) & \( 4.94\) & \(    37\) & \(    47\) & \(    47\) & \(    47\) & \(    57\) & \(    47\) & first\_order \\
JUNKTURN & \( 10010\) & \(  7000\) & \( 9.94\)e\(-06\) & \( 6.79\)e\(-07\) & \( 1.00\)e\(-07\) & \( 5.29\) & \(   108\) & \(   131\) & \(   119\) & \(   131\) & \(   129\) & \(   118\) & first\_order \\
LUKVLE11 & \(  9998\) & \(  6664\) & \( 9.32\)e\(+02\) & \( 4.57\)e\(-08\) & \( 1.00\)e\(-06\) & \( 1.86\) & \(    37\) & \(    63\) & \(    47\) & \(    63\) & \(    57\) & \(    47\) & first\_order \\
LUKVLE17 & \(  9997\) & \(  7497\) & \( 3.24\)e\(+04\) & \( 1.59\)e\(-08\) & \( 1.00\)e\(-07\) & \( 2.45\) & \(    60\) & \(    94\) & \(    78\) & \(    94\) & \(    96\) & \(    78\) & first\_order \\
LUKVLE18 & \(  9997\) & \(  7497\) & \( 1.10\)e\(+04\) & \( 2.00\)e\(-12\) & \( 1.00\)e\(-09\) & \( 3.43\) & \(    60\) & \(    81\) & \(    79\) & \(    81\) & \(    98\) & \(    79\) & first\_order \\
ORTHRDS2 & \(  5003\) & \(  2500\) & \( 7.62\)e\(+02\) & \( 9.96\)e\(-08\) & \( 1.00\)e\(-07\) & \( 1.23\) & \(    47\) & \(    63\) & \(    62\) & \(    63\) & \(    77\) & \(    62\) & first\_order
\\ \hline
\et
\hss}
\end{table}


\section{Nonlinear least squares}
\label{sec:nls}

An important class of problems worthy of special attention is \emph{nonlinear least-squares (NLS) problems} of the form
\begin{equation}  \label{eq:NLS}
    \min_x\ \half \norm{c(x)}^2 \text{ \ subject to \ } \ell \le x \le u,
\end{equation}
where the Jacobian of $c(x)$ is again $J(x)$, and
the bounds are often empty.  Such problems are not immediately meaningful
to Algorithm NCL, but if they are presented in the (probably infeasible) form
\begin{equation}  \label{eq:NLS2}
    \min_x\ 0 \text{ \ subject to \ } c(x) = 0, \quad \ell \le x \le u,
\end{equation}
the first NCL subproblem will be
\[
   \Problem{\NCa}{x,\,r}
   {y_0\T r + \half \rho_0 \norm{r}^2}
   {c(x) + r = 0, \quad \ell \le x \le u,}
\]
which is well suited to \KNITRO and is equivalent to \eqref{eq:NLS} if $y_0=0$ and $\rho_0 > 0$.
If we treat NLS problems as a special case, we can set $y_0=0$, $\rho_0=1$, $\eta_0 = \etastar$, $\omega_0 = \omegastar$ and obtain an optimal solution in one \NCL iteration.
In this sense, Algorithm \NCL is ideally suited to NLS problems \eqref{eq:NLS}.


The \CUTEst collection features a number of NLS problems in both forms~\eqref{eq:NLS} and~\eqref{eq:NLS2}.
While formulation~\eqref{eq:NLS} allows evaluation of the objective gradient
$J(x)^T c(x)$, it does not give access to $J(x)$ itself.  In contrast,
a problem modeled as~\eqref{eq:NLS2} allows solvers to access $J(x)$ directly.

The NLPModels modeling package allows us to formulate~\eqref{eq:NLS} from a problem given as~\eqref{eq:NLS2} and fulfill requests for $J(x)$ in~\eqref{eq:NLS} by returning the constraint Jacobian of~\eqref{eq:NLS2}.
Alternatively, problem~\NCa{} is easily created by the \texttt{NCLModel} constructor.
The construction of both models is illustrated in Listing~\ref{lst:feasibility-residual}.
Once a problem in the form~\eqref{eq:NLS} has been simulated in this way, it can be passed to \KNITRO's nonlinear least-squares solver, which is a variant of the Levenberg-Marquardt method in which bound constraints are treated via an interior-point method.

\lstset{escapeinside={(*@}{@*)}}  
\begin{lstlisting}[label={lst:feasibility-residual}, caption={Formulating~\eqref{eq:NLS} from~\eqref{eq:NLS2}}]
julia> using CUTEst
julia> model = CUTEstModel("ARWHDNE")          # problem in the form (*@\eqref{eq:NLS2}@*)
julia> nls_model = FeasibilityResidual(model)  # interpretation of (*@\eqref{eq:NLS2}@*) as representing (*@\eqref{eq:NLS}@*)
julia> knitro(nls_model)  # NLPModelsKnitro calls KNITRO/Levenberg-Marquardt
julia> ncl_model = NCLModel(model, y=zeros(model.meta.ncon), ρ=1.0)  # problem (*@\NCa@*)
julia> knitro(ncl_model)  # NLPModelsKnitro calls standard KNITRO
\end{lstlisting}

We identified \(127\) problems in the form~\eqref{eq:NLS2} in \CUTEst.
We solve each problem in two ways:
\bc
\bt{ll}
   Solver \knitronls & applies \KNITRO's nonlinear least-squares method to~\eqref{eq:NLS}.
\\ Solver \nclnls    & uses \KNITRO to perform a single \NCL iteration on \NCa.
\et
\ec
In both cases, \KNITRO is given a maximum of \(500\) iterations and \(30\) minutes of CPU time.
Optimality and feasibility tolerances are set to \(10^{-6}\).

\knitronls solved \(101\) problems to optimality, reached the iteration limit in \(19\) cases and the time limit in \(3\) cases, and failed for another reason in \(4\) cases.
\nclnls solved \(119\) problems to optimality, reached the iteration limit in \(3\) cases and the time limit in \(3\) cases, and failed for another reason in \(2\) cases.

\begin{figure}[p]
    \centering
    \includegraphics[height=.3\textheight]{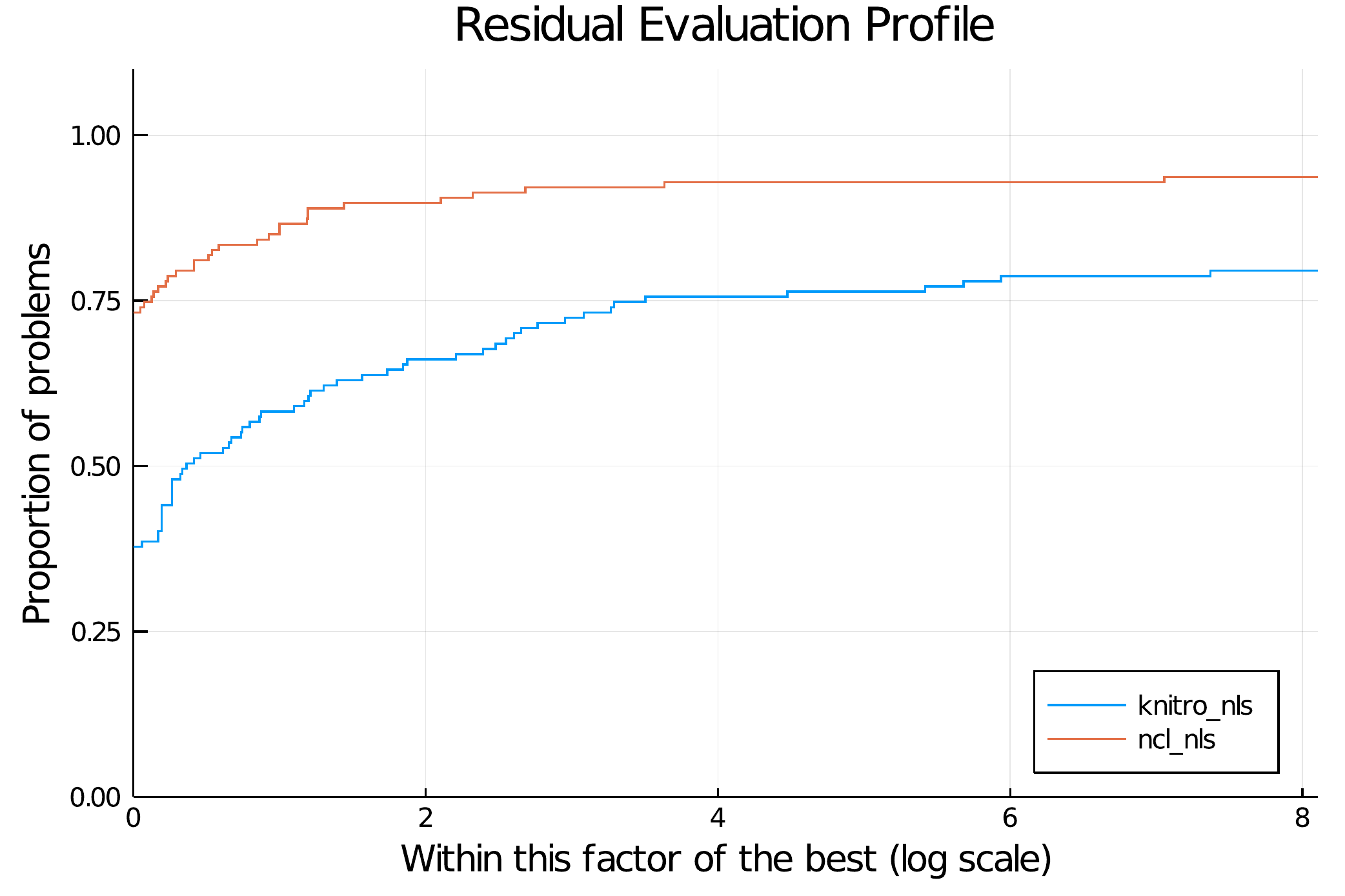}
    \\ \vfill
    \includegraphics[height=.3\textheight]{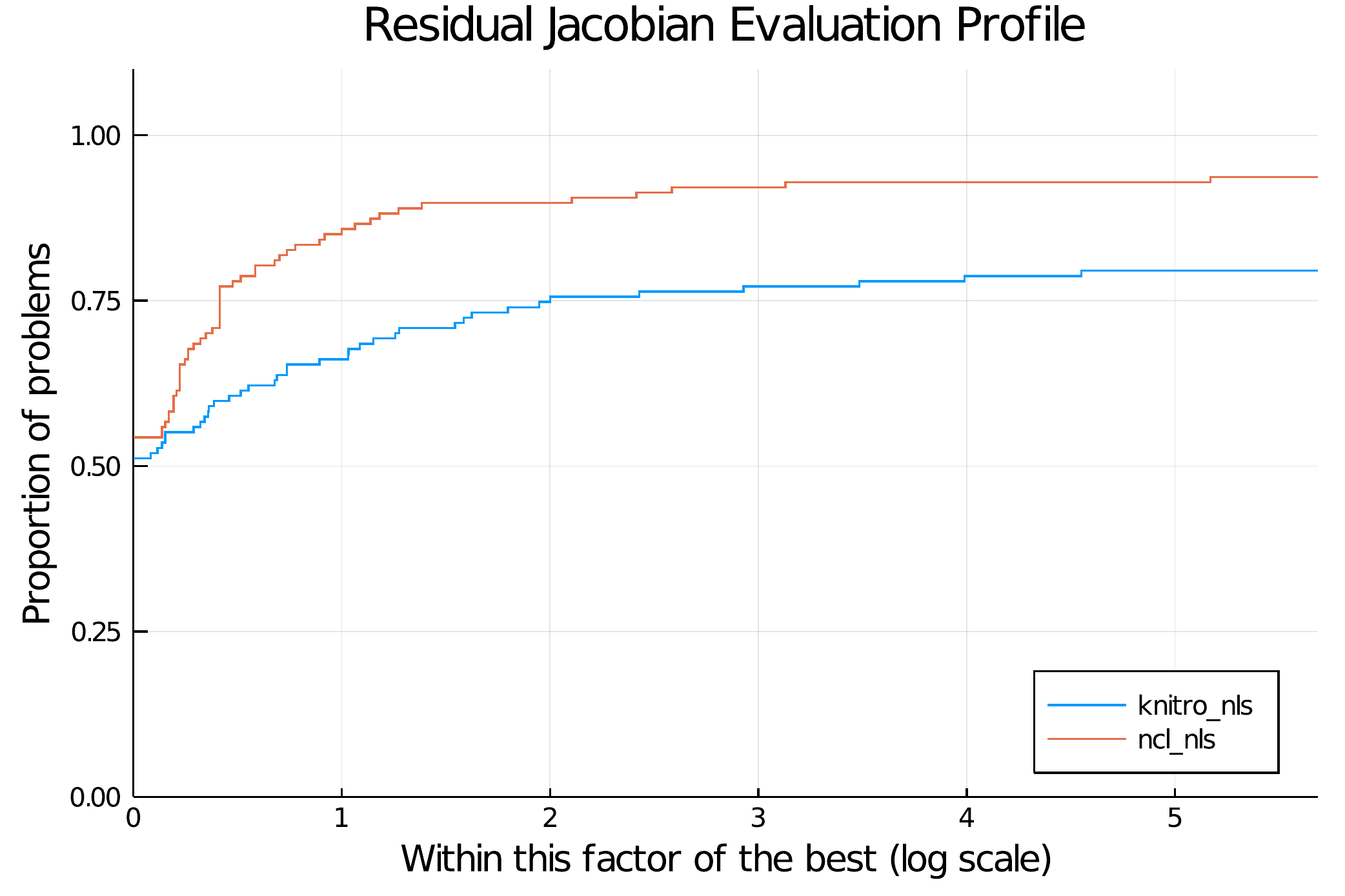}
    \\ \vfill
    \includegraphics[height=.3\textheight]{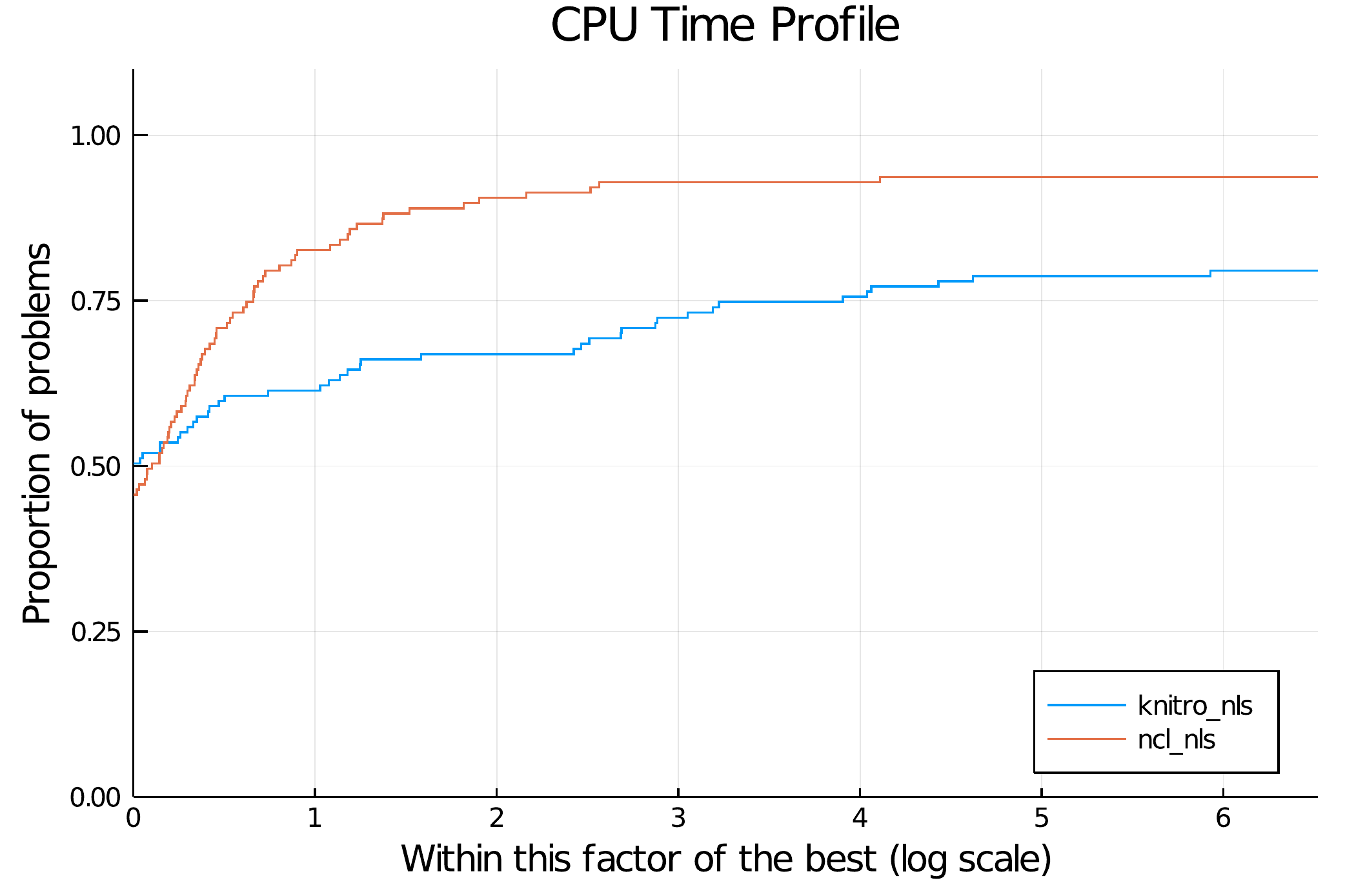}
    \caption{Performance profiles comparing \knitronls
    (KNITRO's NLS solver applied to~\eqref{eq:NLS}) and \nclnls
    (KNITRO solving \NCa)
    on 127 nonlinear least squares problems from CUTEst.
    \nclnls is more efficient.}
    \label{fig:nls-profiles}
\end{figure}

Figure~\ref{fig:nls-profiles} shows Dolan-Mor\'e performance profiles comparing the two solvers.
The top and middle plots use the number of residual and residual Jacobian evaluations as metric, which, in the case of~\eqref{eq:NLS2}, corresponds to the number of constraint and constraint Jacobian evaluations.
The bottom plot uses time as metric.
\nclnls outperforms \knitronls in all three measures and appears substantially more robust.
It is important to keep in mind that a key difference between the two algorithms is that
\nclnls uses second-order information, and therefore performs Hessian evaluations.
Nevertheless, those evaluations are not so costly as to put \NCL at a disadvantage in terms of run-time.  For reference,
Tables~\ref{tab:knitro-nls} and~\ref{tab:ncl-nls} in Appendix~\ref{sec:detailed-results} give the detailed results.

\section{Summary}

Our \AMPL implementation of the tax policy models and Algorithm \NCL has been the only way we could handle these particular problems reliably \cite{NAOIV18}, with \KNITRO solving each subproblem accurately.  Our Julia implementation of \NCL achieves greater efficiency on these \AMPL models by gradually tightening the \KNITRO feasibility and optimality tolerances.
It also permits testing on a broad range of problems, as illustrated on nonlinear least-squares problems and other problems from the \CUTEst test set.
We believe Algorithm \NCL could become an effective general-purpose optimization solver when first and second derivatives are available.  It is especially useful when the LICQ is not satisfied at the solution.
The current Julia implementation of \NCL (with \KNITRO as subproblem solver) is not quite competitive with \KNITRO itself on the general \CUTEst problems in terms of run-time or number of evaluations, but it does solve some problems on which \KNITRO fails.
An advantage is that the implementation is generic and may be applied to problems from any collection adhering to the interface of the NLPModels.jl package~\cite{orban-siqueira-nlpmodels-2020}.


\section*{Acknowledgments}

We are deeply grateful to Professor Ken Judd and Dr Che-Lin Su for
developing the \AMPL tax policy model \cite{JuddSu2011} that led to
the development of Algorithm NCL \cite{NAOIV18}, and to the developers
of \AMPL, Julia, \IPOPT, and \KNITRO for making the implementation
and evaluation of Algorithm NCL possible.  In particular, we thank Dr
Richard Waltz of Artelys for his help in finding runtime options
for warm-starting \KNITRO.  We also give sincere thanks to
Pierre-\'{E}lie Personnaz for obtaining the Julia/NCL results in
Tables \ref{tab:tax1D}--\ref{tab:tax5D}, and to Professor Mehiddin
Al-Baali and other organizers of the NAO-V conference \emph{Numerical
  Analysis and Optimization} at Sultan Qaboos University, Muscat,
Oman, which brought the authors together in January 2020.
Finally, we are very grateful to the referees for their constructive
questions and comments, and to Michael Friedlander for his helpful discussions.

\frenchspacing
\bibliographystyle{plain}
\bibliography{nclNAOV}

\begin{thebibliography}{10}

\bibitem{AMPL}
{AMPL} modeling system.
\newblock \url{http://www.ampl.com}.

\bibitem{bentley-1986}
Jon Bentley.
\newblock Programming pearls: Little languages.
\newblock {\em Commun. ACM}, 29(8):711--–721, August 1986.

\bibitem{bezanson2017julia}
Jeff Bezanson, Alan Edelman, Stefan Karpinski, and Viral~B Shah.
\newblock Julia: A fresh approach to numerical computing.
\newblock {\em SIAM Rev.}, 59(1):65--98, 2017.

\bibitem{CGT91}
A.~R. Conn, N.~I.~M. Gould, and {Ph.}~L. Toint.
\newblock A globally convergent augmented {Lagrangian} algorithm for
  optimization with general constraints and simple bounds.
\newblock {\em SIAM J. Numer. Anal.}, 28:545--572, 1991.

\bibitem{CGT92}
A.~R. Conn, N.~I.~M. Gould, and {Ph.}~L. Toint.
\newblock {\em {LANCELOT: A Fortran Package for Large-scale Nonlinear
  Optimization (Release A)}}.
\newblock Lecture Notes in Computation Mathematics 17. Springer Verlag, Berlin,
  Heidelberg, New York, London, Paris and Tokyo, 1992.

\bibitem{AMPLbook}
R.~Fourer, D.~M. Gay, and B.~W. Kernighan.
\newblock {\em {AMPL}: A Modeling Language for Mathematical Programming}.
\newblock Brooks/Cole, Pacific Grove, second edition, 2002.

\bibitem{GAMS}
{GAMS} modeling system.
\newblock \url{http://www.gams.com}.

\bibitem{GilMS05}
P.~E. Gill, W.~Murray, and M.~A. Saunders.
\newblock {SNOPT}: An {SQP} algorithm for large-scale constrained optimization.
\newblock {\em SIAM Rev.}, 47(1):99--131, 2005.
\newblock SIGEST article.

\bibitem{CUTEst}
N.~I.~M. Gould, D.~Orban, and {Ph.}~L. Toint.
\newblock {CUTEst}: a {Constrained and Unconstrained Testing Environment with
  safe threads}.
\newblock {\em Comput. Optim. Appl.}, 60:545--557, 2015.

\bibitem{IPOPT}
{COIN-OR} {I}nterior {P}oint {O}ptimizer {IPOPT}.
\newblock \url{https://github.com/coin-or/Ipopt}.

\bibitem{JuddSu2011}
K.~L. Judd and C.-L. Su.
\newblock Optimal income taxation with multidimensional taxpayer types.
\newblock Working paper, Hoover Institution, Stanford University, 2011.

\bibitem{KNITRO}
{KNITRO} optimization software.
\newblock \url{https://www.artelys.com/tools/knitro_doc/2_userGuide.html}.

\bibitem{LANCELOT}
{LANCELOT} optimization software.
\newblock \url{http://www.numerical.rl.ac.uk/lancelot/blurb.html}.

\bibitem{jump}
M.~Lubin and I.~Dunning.
\newblock Computing in operations research using julia.
\newblock {\em INFORMS J. Comput.}, 27(2), 2015.

\bibitem{NAOIV18}
D.~Ma, K.~L. Judd, D.~Orban, and M.~A. Saunders.
\newblock Stabilized optimization via an {NCL} algorithm.
\newblock In M.~{Al-Baali et al.}, editor, {\em Numerical Analysis and
  Optimization, NAO-IV, Muscat, Oman, January 2017}, pages 173--191. Springer
  International Publishing AG, 2018.

\bibitem{MurS82}
B.~A. Murtagh and M.~A. Saunders.
\newblock A projected {Lagrangian} algorithm and its implementation for sparse
  nonlinear constraints.
\newblock {\em Math. Program. Study}, 16:84--117, 1982.

\bibitem{NCLweb}
{NCL AMPL models}.
\newblock \url{http://stanford.edu/group/SOL/multiscale/models/NCL/}.

\bibitem{NEOS}
{NEOS server for optimization}.
\newblock \url{http://www.neos-server.org/neos/}.

\bibitem{orban-personnaz-ncl-2020}
D.~Orban and P.~E. Personnaz.
\newblock {NCL.jl}: A nonlinearly-constrained augmented-{L}agrangian method.
\newblock \url{https://github.com/JuliaSmoothOptimizers/NCL.jl}, July 2020.

\bibitem{orban-siqueira-amplnlreader-2020}
D.~Orban, A.~S. Siqueira, and {contributors}.
\newblock {AmplNLReader.jl}: A {J}ulia interface to {AMPL}.
\newblock \url{https://github.com/JuliaSmoothOptimizers/AmplNLReader.jl}, July
  2020.

\bibitem{orban-siqueira-cutest-2020}
D.~Orban, A.~S. Siqueira, and {contributors}.
\newblock {CUTEst.jl}: {J}ulia's {CUTEst} interface.
\newblock \url{https://github.com/JuliaSmoothOptimizers/CUTEst.jl}, October
  2020.

\bibitem{orban-siqueira-jso-2020}
D.~Orban, A.~S. Siqueira, and {contributors}.
\newblock {JuliaSmoothOptimizers}: Infrastructure and solvers for continuous
  optimization in {J}ulia.
\newblock \url{https://github.com/JuliaSmoothOptimizers}, July 2020.

\bibitem{orban-siqueira-nlpmodelsipopt-2020}
D.~Orban, A.~S. Siqueira, and {contributors}.
\newblock {NLPModelsIpopt.jl}: A thin {IPOPT} wrapper for {NLPModels}.
\newblock \url{https://github.com/JuliaSmoothOptimizers/NLPModelsIpopt.jl},
  July 2020.

\bibitem{orban-siqueira-nlpmodels-2020}
D.~Orban, A.~S. Siqueira, and {contributors}.
\newblock {NLPModels.jl}: Data structures for optimization models.
\newblock \url{https://github.com/JuliaSmoothOptimizers/NLPModels.jl}, July
  2020.

\bibitem{orban-siqueira-nlpmodelsknitro-2020}
D.~Orban, A.~S. Siqueira, and {contributors}.
\newblock {NLPModelsKnitro.jl}: A thin {KNITRO} wrapper for {NLPModels}.
\newblock \url{https://github.com/JuliaSmoothOptimizers/NLPModelsKnitro.jl},
  July 2020.

\end{thebibliography}

\appendix
\section{Detailed results for Julia/NCL on NLS problems}
\label{sec:detailed-results}

Table~\ref{tab:knitro-nls} reports the detailed results of KNITRO/Levenberg-Marquardt on problems of the form~\eqref{eq:NLS} using the modeling mechanism of Section~\ref{sec:nls}.
In the table headers, ``nvar'' is the number of variables, ``ncon'' is the number of constraints (i.e., the number of least-squares residuals), \(f\) is the final objective value, \(\|\nabla L\|_2\) is the final dual residual, \(t\) is the run-time in seconds, ``iter'' is the number of iterations, ``\#\(c\)'' is the number of constraint (i.e, residual) evaluations, ``\#\(\nabla c\)'' is the number of constraint (i.e., residual) Jacobian evaluations, and ``status'' is the final solver status.

Table~\ref{tab:ncl-nls} reports the results of Julia/NCL solving Problem~\NCa for the same models.
In the interest of space, the second table does not repeat problem dimensions.
The other columns are as follows: \(\|c\|_2\) is the final primal feasibility, and \(\#\nabla^2 L\) is the number of Hessian evaluations.

\begin{longtable}{lrrrrrrrrl}
\caption{%
    \label{tab:knitro-nls}
    \knitronls results on \(127\) CUTEst nonlinear least-squares problems.
    \(101\) problems were solved successfully.
}\\
\hline
name & nvar & ncon & \(f\)\hspace*{14pt} & \(\|\nabla L\|_2\) & \(t\)\hspace*{6pt} & iter & \#\(c\) & \#\(\nabla c\) & status \\\hline
\endhead
\hline
\multicolumn{10}{r}{{\bfseries Continued on next page}}\\
\hline
\endfoot
\endlastfoot
ARWHDNE & \(   500\) & \(   998\) & \( 6.971\)e\(+01\) & \(8.5\)e\(-06\) & \(0.46\) & \(    22\) & \(   167\) & \(    23\) & first\_order \\
BA-L1 & \(    57\) & \(    12\) & \( 1.204\)e\(-25\) & \(3.8\)e\(-10\) & \(0.00\) & \(     5\) & \(     6\) & \(     6\) & first\_order \\
BA-L16 & \( 66462\) & \(167436\) & \( 4.243\)e\(+05\) & \(4.2\)e\(-04\) & \(69.92\) & \(    27\) & \(    30\) & \(    28\) & first\_order \\
BA-L1SP & \(    57\) & \(    12\) & \( 1.620\)e\(-23\) & \(5.1\)e\(-09\) & \(0.01\) & \(     5\) & \(     6\) & \(     6\) & first\_order \\
BA-L21 & \( 34134\) & \( 72910\) & \( 1.975\)e\(+05\) & \(7.2\)e\(-05\) & \(129.28\) & \(   119\) & \(   160\) & \(   120\) & first\_order \\
BA-L49 & \( 23769\) & \( 63686\) & \( 1.241\)e\(+05\) & \(4.0\)e\(+05\) & \(1809.44\) & \(   121\) & \(   728\) & \(   122\) & max\_time \\
BA-L52 & \(192627\) & \(694346\) & \( 2.147\)e\(+06\) & \(3.9\)e\(+07\) & \(1808.91\) & \(    37\) & \(   200\) & \(    38\) & max\_time \\
BA-L73 & \( 33753\) & \( 92244\) & \( 6.312\)e\(+05\) & \(9.5\)e\(+06\) & \(1801.24\) & \(   396\) & \(  2312\) & \(   397\) & max\_time \\
BARDNE & \(     3\) & \(    15\) & \( 4.107\)e\(-03\) & \(2.7\)e\(-09\) & \(0.00\) & \(     5\) & \(     6\) & \(     6\) & first\_order \\
BDQRTICNE & \(  5000\) & \(  9992\) & \( 1.000\)e\(+04\) & \(1.4\)e\(-01\) & \(91.67\) & \(   500\) & \(  3668\) & \(   501\) & max\_iter \\
BEALENE & \(     2\) & \(     3\) & \( 1.491\)e\(-25\) & \(2.6\)e\(-12\) & \(0.00\) & \(     6\) & \(     8\) & \(     7\) & first\_order \\
BIGGS6NE & \(     6\) & \(    13\) & \( 2.367\)e\(-17\) & \(5.9\)e\(-09\) & \(0.01\) & \(    30\) & \(   124\) & \(    31\) & first\_order \\
BOX3NE & \(     3\) & \(    10\) & \( 3.259\)e\(-19\) & \(3.9\)e\(-10\) & \(0.00\) & \(     5\) & \(     6\) & \(     6\) & first\_order \\
BROWNBSNE & \(     2\) & \(     3\) & \( 0.000\)e\(+00\) & \(0.0\)e\(+00\) & \(0.00\) & \(    12\) & \(    53\) & \(    13\) & first\_order \\
BROWNDENE & \(     4\) & \(    20\) & \( 4.291\)e\(+04\) & \(2.4\)e\(+00\) & \(0.09\) & \(   500\) & \(  1697\) & \(   501\) & max\_iter \\
BRYBNDNE & \(  5000\) & \(  5000\) & \( 1.499\)e\(-21\) & \(8.2\)e\(-11\) & \(1.16\) & \(     6\) & \(     7\) & \(     7\) & first\_order \\
CHAINWOONE & \(  4000\) & \( 11994\) & \( 4.699\)e\(+03\) & \(5.9\)e\(+01\) & \(63.96\) & \(   500\) & \(  1339\) & \(   501\) & max\_iter \\
CHEBYQADNE & \(   100\) & \(   100\) & \( 4.749\)e\(-03\) & \(1.5\)e\(-04\) & \(7.53\) & \(   500\) & \(  1994\) & \(   501\) & max\_iter \\
CHNRSBNE & \(    50\) & \(    98\) & \( 7.394\)e\(-18\) & \(3.2\)e\(-08\) & \(0.01\) & \(    38\) & \(    78\) & \(    39\) & first\_order \\
CHNRSNBMNE & \(    50\) & \(    98\) & \( 1.452\)e\(-20\) & \(1.8\)e\(-09\) & \(0.02\) & \(    54\) & \(   132\) & \(    55\) & first\_order \\
COATINGNE & \(   134\) & \(   252\) & \( 2.527\)e\(-01\) & \(5.5\)e\(-07\) & \(0.01\) & \(     9\) & \(    11\) & \(    10\) & first\_order \\
CUBENE & \(     2\) & \(     2\) & \( 1.085\)e\(-26\) & \(1.1\)e\(-13\) & \(0.00\) & \(     4\) & \(     9\) & \(     5\) & first\_order \\
DECONVBNE & \(    63\) & \(    40\) & \( 4.542\)e\(-10\) & \(6.0\)e\(-07\) & \(0.03\) & \(    54\) & \(   214\) & \(    55\) & first\_order \\
DECONVNE & \(    63\) & \(    40\) & \( 2.245\)e\(-16\) & \(8.5\)e\(-10\) & \(0.00\) & \(     2\) & \(     3\) & \(     3\) & first\_order \\
DENSCHNBNE & \(     2\) & \(     3\) & \( 8.573\)e\(-27\) & \(1.9\)e\(-13\) & \(0.00\) & \(     5\) & \(     6\) & \(     6\) & first\_order \\
DENSCHNCNE & \(     2\) & \(     2\) & \( 2.838\)e\(-22\) & \(8.4\)e\(-11\) & \(0.00\) & \(     7\) & \(     8\) & \(     8\) & first\_order \\
DENSCHNDNE & \(     3\) & \(     3\) & \( 8.307\)e\(-10\) & \(5.5\)e\(-07\) & \(0.00\) & \(    17\) & \(    18\) & \(    18\) & first\_order \\
DENSCHNENE & \(     3\) & \(     3\) & \( 3.005\)e\(-22\) & \(2.4\)e\(-11\) & \(0.00\) & \(     8\) & \(    19\) & \(     9\) & first\_order \\
DENSCHNFNE & \(     2\) & \(     2\) & \( 1.392\)e\(-26\) & \(2.1\)e\(-12\) & \(0.00\) & \(     5\) & \(     6\) & \(     6\) & first\_order \\
DEVGLA1NE & \(     4\) & \(    24\) & \( 1.063\)e\(-13\) & \(2.4\)e\(-08\) & \(0.00\) & \(    11\) & \(    31\) & \(    12\) & first\_order \\
DEVGLA2NE & \(     5\) & \(    16\) & \( 4.828\)e\(-15\) & \(5.9\)e\(-07\) & \(0.00\) & \(     9\) & \(    16\) & \(    10\) & first\_order \\
EGGCRATENE & \(     2\) & \(     4\) & \( 4.744\)e\(+00\) & \(9.6\)e\(-07\) & \(0.00\) & \(     5\) & \(     6\) & \(     6\) & first\_order \\
ELATVIDUNE & \(     2\) & \(     3\) & \( 2.738\)e\(+01\) & \(3.2\)e\(-06\) & \(0.00\) & \(    15\) & \(    30\) & \(    16\) & first\_order \\
ENGVAL2NE & \(     3\) & \(     5\) & \( 2.465\)e\(-32\) & \(2.2\)e\(-16\) & \(0.00\) & \(    10\) & \(    14\) & \(    11\) & first\_order \\
ERRINROSNE & \(    50\) & \(    98\) & \( 2.020\)e\(+01\) & \(1.8\)e\(-05\) & \(0.01\) & \(    45\) & \(    63\) & \(    46\) & first\_order \\
ERRINRSMNE & \(    50\) & \(    98\) & \( 1.926\)e\(+01\) & \(1.7\)e\(-05\) & \(0.01\) & \(    44\) & \(    66\) & \(    45\) & first\_order \\
EXP2NE & \(     2\) & \(    10\) & \( 1.618\)e\(-19\) & \(8.7\)e\(-12\) & \(0.00\) & \(     5\) & \(     6\) & \(     6\) & first\_order \\
EXPFITNE & \(     2\) & \(    10\) & \( 1.203\)e\(-01\) & \(4.9\)e\(-07\) & \(0.00\) & \(    10\) & \(    12\) & \(    11\) & first\_order \\
EXTROSNBNE & \(  1000\) & \(   999\) & \( 8.071\)e\(-31\) & \(1.3\)e\(-14\) & \(0.07\) & \(     6\) & \(    14\) & \(     7\) & first\_order \\
FBRAIN2NE & \(     4\) & \(  2211\) & \( 1.842\)e\(-01\) & \(6.8\)e\(-07\) & \(0.08\) & \(     8\) & \(    11\) & \(     9\) & first\_order \\
FBRAINNE & \(     2\) & \(  2211\) & \( 2.083\)e\(-01\) & \(5.7\)e\(-07\) & \(0.02\) & \(     5\) & \(     6\) & \(     6\) & first\_order \\
FREURONE & \(     2\) & \(     2\) & \( 2.449\)e\(+01\) & \(1.8\)e\(-05\) & \(0.01\) & \(    19\) & \(   102\) & \(    20\) & first\_order \\
GENROSEBNE & \(   500\) & \(   998\) & \( 7.965\)e\(+02\) & \(3.6\)e\(-06\) & \(0.04\) & \(     8\) & \(    10\) & \(    10\) & first\_order \\
GENROSENE & \(  1000\) & \(  1999\) & \( 1.247\)e\(+02\) & \(6.3\)e\(+00\) & \(3.40\) & \(   500\) & \(  1567\) & \(   501\) & max\_iter \\
GULFNE & \(     3\) & \(    99\) & \( 2.107\)e\(-01\) & \(4.0\)e\(+06\) & \(0.27\) & \(   500\) & \(  2578\) & \(   501\) & max\_iter \\
HATFLDANE & \(     4\) & \(     4\) & \( 2.172\)e\(-17\) & \(3.4\)e\(-09\) & \(0.00\) & \(     8\) & \(    10\) & \(     9\) & first\_order \\
HATFLDBNE & \(     4\) & \(     4\) & \( 2.786\)e\(-03\) & \(1.4\)e\(-07\) & \(0.00\) & \(     7\) & \(     8\) & \(     8\) & first\_order \\
HATFLDCNE & \(    25\) & \(    25\) & \( 1.018\)e\(-15\) & \(1.5\)e\(-08\) & \(0.00\) & \(     3\) & \(     4\) & \(     4\) & first\_order \\
HATFLDDNE & \(     3\) & \(    10\) & \( 1.273\)e\(-07\) & \(7.1\)e\(-11\) & \(0.00\) & \(     6\) & \(     9\) & \(     7\) & first\_order \\
HATFLDENE & \(     3\) & \(    21\) & \( 1.364\)e\(-06\) & \(8.2\)e\(-12\) & \(0.00\) & \(     5\) & \(     6\) & \(     6\) & first\_order \\
HATFLDFLNE & \(     3\) & \(     3\) & \( 3.254\)e\(-05\) & \(9.1\)e\(-07\) & \(0.01\) & \(    11\) & \(    69\) & \(    12\) & first\_order \\
HELIXNE & \(     3\) & \(     3\) & \( 2.195\)e\(-20\) & \(2.8\)e\(-09\) & \(0.00\) & \(     9\) & \(    11\) & \(    10\) & first\_order \\
HIMMELBFNE & \(     4\) & \(     7\) & \( 1.593\)e\(+06\) & \(4.4\)e\(-04\) & \(0.01\) & \(    28\) & \(    66\) & \(    29\) & first\_order \\
HS1NE & \(     2\) & \(     2\) & \( 3.274\)e\(-17\) & \(4.0\)e\(-09\) & \(0.00\) & \(     9\) & \(    20\) & \(    10\) & first\_order \\
HS25NE & \(     3\) & \(    99\) & \( 1.642\)e\(+01\) & \(9.2\)e\(-09\) & \(0.00\) & \(     0\) & \(     1\) & \(     1\) & first\_order \\
HS2NE & \(     2\) & \(     2\) & \( 2.471\)e\(+00\) & \(2.1\)e\(-08\) & \(0.00\) & \(     6\) & \(     8\) & \(     8\) & first\_order \\
INTEQNE & \(    12\) & \(    12\) & \( 6.938\)e\(-14\) & \(2.1\)e\(-07\) & \(0.00\) & \(     2\) & \(     3\) & \(     3\) & first\_order \\
JENSMPNE & \(     2\) & \(    10\) & \( 6.218\)e\(+01\) & \(7.1\)e\(-06\) & \(0.01\) & \(    10\) & \(    82\) & \(    11\) & first\_order \\
JUDGENE & \(     2\) & \(    20\) & \( 8.041\)e\(+00\) & \(1.2\)e\(-06\) & \(0.00\) & \(     9\) & \(    10\) & \(    10\) & first\_order \\
KOEBHELBNE & \(     3\) & \(   156\) & \( 3.876\)e\(+01\) & \(7.5\)e\(-07\) & \(0.10\) & \(   194\) & \(   757\) & \(   195\) & first\_order \\
KOWOSBNE & \(     4\) & \(    11\) & \( 1.539\)e\(-04\) & \(6.9\)e\(-07\) & \(0.00\) & \(    16\) & \(    30\) & \(    17\) & first\_order \\
LIARWHDNE & \(  5000\) & \( 10000\) & \( 1.608\)e\(-27\) & \(4.4\)e\(-12\) & \(1.07\) & \(     5\) & \(     6\) & \(     6\) & first\_order \\
LINVERSENE & \(  1999\) & \(  2997\) & \( 3.405\)e\(+02\) & \(1.2\)e\(-02\) & \(13.28\) & \(   500\) & \(  1753\) & \(   502\) & max\_iter \\
MANCINONE & \(   100\) & \(   100\) & \( 7.966\)e\(-22\) & \(2.2\)e\(-08\) & \(0.10\) & \(     5\) & \(     6\) & \(     6\) & first\_order \\
MANNE & \(  6000\) & \(  4000\) & \( 3.261\)e\(+39\) & \(1.4\)e\(+19\) & \(95.71\) & \(   167\) & \(   168\) & \(   168\) & unknown \\
MARINE & \( 11215\) & \( 11192\) & \( 3.922\)e\(-21\) & \(8.0\)e\(-08\) & \(8.51\) & \(     9\) & \(    10\) & \(    10\) & first\_order \\
MEYER3NE & \(     3\) & \(    16\) & \( 4.397\)e\(+01\) & \(1.1\)e\(-03\) & \(0.00\) & \(    11\) & \(    19\) & \(    12\) & unknown \\
MODBEALENE & \( 20000\) & \( 39999\) & \( 1.374\)e\(+00\) & \(8.3\)e\(-07\) & \(109.80\) & \(    38\) & \(    73\) & \(    39\) & first\_order \\
MOREBVNE & \(    10\) & \(    10\) & \( 1.085\)e\(-14\) & \(1.9\)e\(-08\) & \(0.00\) & \(     2\) & \(     3\) & \(     3\) & first\_order \\
MUONSINE & \(     1\) & \(   512\) & \( 2.194\)e\(+04\) & \(8.8\)e\(-04\) & \(0.02\) & \(    36\) & \(    37\) & \(    37\) & first\_order \\
NGONE & \(   200\) & \(  5048\) & \( 7.203\)e\(-13\) & \(8.3\)e\(-07\) & \(36.90\) & \(   221\) & \(   914\) & \(   223\) & first\_order \\
NONDIANE & \(  5000\) & \(  5000\) & \( 4.949\)e\(-01\) & \(1.9\)e\(-08\) & \(2.63\) & \(    16\) & \(    39\) & \(    17\) & first\_order \\
NONMSQRTNE & \(  4900\) & \(  4900\) & \( 3.543\)e\(+02\) & \(1.1\)e\(+01\) & \(154.67\) & \(   500\) & \(  2624\) & \(   501\) & max\_iter \\
NONSCOMPNE & \(  5000\) & \(  5000\) & \( 1.378\)e\(-06\) & \(5.6\)e\(-07\) & \(11.31\) & \(    80\) & \(   217\) & \(    81\) & first\_order \\
OSCIGRNE & \(100000\) & \(100000\) & \( 3.138\)e\(-24\) & \(8.3\)e\(-09\) & \(3.62\) & \(     7\) & \(     8\) & \(     8\) & first\_order \\
OSCIPANE & \(    10\) & \(    10\) & \( 5.000\)e\(-01\) & \(1.1\)e\(-01\) & \(0.11\) & \(   500\) & \(  2880\) & \(   501\) & max\_iter \\
PALMER1ANE & \(     6\) & \(    35\) & \( 4.494\)e\(-02\) & \(6.7\)e\(-07\) & \(0.01\) & \(    25\) & \(    84\) & \(    26\) & first\_order \\
PALMER1BNE & \(     4\) & \(    35\) & \( 1.724\)e\(+00\) & \(2.6\)e\(-08\) & \(0.00\) & \(     6\) & \(     7\) & \(     7\) & first\_order \\
PALMER1ENE & \(     8\) & \(    35\) & \( 1.822\)e\(-01\) & \(1.3\)e\(+02\) & \(0.17\) & \(   500\) & \(  3471\) & \(   501\) & max\_iter \\
PALMER1NE & \(     4\) & \(    31\) & \( 5.877\)e\(+03\) & \(8.4\)e\(+00\) & \(0.13\) & \(   500\) & \(  2561\) & \(   501\) & max\_iter \\
PALMER2ANE & \(     6\) & \(    23\) & \( 8.555\)e\(-03\) & \(5.9\)e\(-07\) & \(0.02\) & \(    58\) & \(   248\) & \(    59\) & first\_order \\
PALMER2BNE & \(     4\) & \(    23\) & \( 3.116\)e\(-01\) & \(5.2\)e\(-08\) & \(0.00\) & \(     8\) & \(    10\) & \(     9\) & first\_order \\
PALMER2ENE & \(     8\) & \(    23\) & \( 1.952\)e\(-02\) & \(8.9\)e\(-07\) & \(0.01\) & \(     9\) & \(    68\) & \(    10\) & first\_order \\
PALMER2NE & \(     4\) & \(    23\) & \( 1.826\)e\(+03\) & \(4.1\)e\(-06\) & \(0.01\) & \(    26\) & \(    92\) & \(    27\) & unknown \\
PALMER3ANE & \(     6\) & \(    23\) & \( 1.022\)e\(-02\) & \(1.0\)e\(-07\) & \(0.00\) & \(    18\) & \(    39\) & \(    19\) & first\_order \\
PALMER3BNE & \(     4\) & \(    23\) & \( 2.114\)e\(+00\) & \(1.5\)e\(-07\) & \(0.00\) & \(    11\) & \(    14\) & \(    12\) & first\_order \\
PALMER3ENE & \(     8\) & \(    23\) & \( 2.303\)e\(-02\) & \(1.2\)e\(+01\) & \(0.15\) & \(   500\) & \(  3348\) & \(   501\) & max\_iter \\
PALMER3NE & \(     4\) & \(    23\) & \( 1.133\)e\(+03\) & \(5.7\)e\(-02\) & \(0.15\) & \(   500\) & \(  3303\) & \(   501\) & max\_iter \\
PALMER4ANE & \(     6\) & \(    23\) & \( 2.030\)e\(-02\) & \(3.8\)e\(-07\) & \(0.03\) & \(    98\) & \(   445\) & \(    99\) & first\_order \\
PALMER4BNE & \(     4\) & \(    23\) & \( 3.418\)e\(+00\) & \(6.1\)e\(-07\) & \(0.00\) & \(    14\) & \(    18\) & \(    15\) & first\_order \\
PALMER4ENE & \(     8\) & \(    23\) & \( 6.286\)e\(-02\) & \(9.3\)e\(-07\) & \(0.04\) & \(    69\) & \(   667\) & \(    70\) & first\_order \\
PALMER4NE & \(     4\) & \(    23\) & \( 1.143\)e\(+03\) & \(2.2\)e\(-05\) & \(0.01\) & \(    35\) & \(   163\) & \(    36\) & unknown \\
PALMER5ANE & \(     8\) & \(    12\) & \( 1.706\)e\(-01\) & \(1.2\)e\(+00\) & \(0.13\) & \(   500\) & \(  3000\) & \(   501\) & max\_iter \\
PALMER5BNE & \(     9\) & \(    12\) & \( 4.876\)e\(-03\) & \(3.4\)e\(-07\) & \(0.03\) & \(    93\) & \(   481\) & \(    94\) & first\_order \\
PALMER5ENE & \(     8\) & \(    12\) & \( 1.699\)e\(-02\) & \(1.4\)e\(+00\) & \(0.16\) & \(   500\) & \(  3508\) & \(   501\) & max\_iter \\
PALMER6ANE & \(     6\) & \(    13\) & \( 2.797\)e\(-02\) & \(2.3\)e\(-07\) & \(0.00\) & \(    17\) & \(    30\) & \(    18\) & first\_order \\
PALMER6ENE & \(     8\) & \(    13\) & \( 2.424\)e\(-02\) & \(5.8\)e\(-07\) & \(0.08\) & \(   142\) & \(  1654\) & \(   143\) & first\_order \\
PALMER7ANE & \(     6\) & \(    13\) & \( 5.526\)e\(+00\) & \(2.4\)e\(+00\) & \(0.13\) & \(   500\) & \(  2806\) & \(   501\) & max\_iter \\
PALMER7ENE & \(     8\) & \(    13\) & \( 3.353\)e\(+00\) & \(1.4\)e\(+02\) & \(0.16\) & \(   500\) & \(  3243\) & \(   500\) & max\_iter \\
PALMER8ANE & \(     6\) & \(    12\) & \( 3.700\)e\(-02\) & \(5.8\)e\(-07\) & \(0.02\) & \(    72\) & \(   303\) & \(    73\) & first\_order \\
PALMER8ENE & \(     8\) & \(    12\) & \( 1.687\)e\(-01\) & \(2.3\)e\(-04\) & \(0.22\) & \(   500\) & \(  5007\) & \(   501\) & max\_iter \\
PENLT1NE & \(    10\) & \(    11\) & \( 3.576\)e\(-10\) & \(8.1\)e\(-07\) & \(0.02\) & \(   111\) & \(   385\) & \(   112\) & first\_order \\
PENLT2NE & \(     4\) & \(     8\) & \( 4.701\)e\(-11\) & \(9.6\)e\(-07\) & \(0.03\) & \(   163\) & \(   613\) & \(   164\) & first\_order \\
PINENE & \(  8805\) & \(  8795\) & \( 5.184\)e\(-17\) & \(1.1\)e\(-07\) & \(1.72\) & \(     2\) & \(    10\) & \(     3\) & first\_order \\
POWERSUMNE & \(     4\) & \(     4\) & \( 2.325\)e\(-17\) & \(7.6\)e\(-07\) & \(0.00\) & \(    17\) & \(    21\) & \(    18\) & first\_order \\
PRICE3NE & \(     2\) & \(     2\) & \( 5.614\)e\(-22\) & \(3.9\)e\(-10\) & \(0.00\) & \(     7\) & \(     8\) & \(     8\) & first\_order \\
PRICE4NE & \(     2\) & \(     2\) & \( 1.321\)e\(-14\) & \(9.8\)e\(-07\) & \(0.00\) & \(    21\) & \(    22\) & \(    22\) & first\_order \\
QINGNE & \(   100\) & \(   100\) & \( 9.461\)e\(-20\) & \(1.2\)e\(-09\) & \(0.00\) & \(     5\) & \(     8\) & \(     6\) & first\_order \\
RSNBRNE & \(     2\) & \(     2\) & \( 4.832\)e\(-30\) & \(3.1\)e\(-15\) & \(0.00\) & \(    11\) & \(    34\) & \(    12\) & first\_order \\
S308NE & \(     2\) & \(     3\) & \( 3.866\)e\(-01\) & \(9.1\)e\(-07\) & \(0.00\) & \(    34\) & \(    36\) & \(    35\) & first\_order \\
SBRYBNDNE & \(  5000\) & \(  5000\) & \( 1.790\)e\(-21\) & \(3.0\)e\(-07\) & \(1.16\) & \(     6\) & \(     7\) & \(     7\) & first\_order \\
SINVALNE & \(     2\) & \(     2\) & \( 3.852\)e\(-32\) & \(2.8\)e\(-15\) & \(0.00\) & \(     4\) & \(    11\) & \(     5\) & first\_order \\
SPECANNE & \(     9\) & \( 15000\) & \( 3.291\)e\(-13\) & \(5.3\)e\(-08\) & \(0.08\) & \(     6\) & \(     7\) & \(     7\) & first\_order \\
SROSENBRNE & \(  5000\) & \(  5000\) & \( 5.547\)e\(-28\) & \(6.7\)e\(-16\) & \(0.49\) & \(     2\) & \(     3\) & \(     3\) & first\_order \\
SSBRYBNDNE & \(  5000\) & \(  5000\) & \( 5.274\)e\(-22\) & \(4.0\)e\(-10\) & \(1.15\) & \(     6\) & \(     7\) & \(     7\) & first\_order \\
STREGNE & \(     4\) & \(     2\) & \( 2.236\)e\(-03\) & \(4.7\)e\(-01\) & \(0.05\) & \(   500\) & \(   535\) & \(   501\) & max\_iter \\
STRTCHDVNE & \(    10\) & \(     9\) & \( 3.723\)e\(-10\) & \(4.4\)e\(-07\) & \(0.00\) & \(     8\) & \(     9\) & \(     9\) & first\_order \\
TQUARTICNE & \(  5000\) & \(  5000\) & \( 2.631\)e\(-26\) & \(2.3\)e\(-13\) & \(0.39\) & \(     1\) & \(     2\) & \(     2\) & first\_order \\
TRIGON1NE & \(    10\) & \(    10\) & \( 6.724\)e\(-16\) & \(1.0\)e\(-07\) & \(0.00\) & \(     4\) & \(     5\) & \(     5\) & first\_order \\
TRIGON2NE & \(    10\) & \(    31\) & \( 1.620\)e\(+00\) & \(1.5\)e\(-09\) & \(0.00\) & \(     7\) & \(     8\) & \(     8\) & first\_order \\
VARDIMNE & \(    10\) & \(    12\) & \( 7.960\)e\(-16\) & \(4.0\)e\(-07\) & \(0.00\) & \(     9\) & \(    10\) & \(    10\) & first\_order \\
VIBRBEAMNE & \(     8\) & \(    30\) & \( 7.822\)e\(-02\) & \(6.9\)e\(-07\) & \(0.00\) & \(    10\) & \(    11\) & \(    11\) & first\_order \\
WATSONNE & \(    12\) & \(    31\) & \( 1.429\)e\(-15\) & \(1.3\)e\(-13\) & \(0.00\) & \(     4\) & \(     5\) & \(     5\) & first\_order \\
WAYSEA1NE & \(     2\) & \(     2\) & \( 3.683\)e\(-16\) & \(8.7\)e\(-07\) & \(0.00\) & \(     7\) & \(     8\) & \(     8\) & first\_order \\
WAYSEA2NE & \(     2\) & \(     2\) & \( 3.388\)e\(-18\) & \(1.3\)e\(-08\) & \(0.00\) & \(    11\) & \(    17\) & \(    12\) & first\_order \\
WEEDSNE & \(     3\) & \(    12\) & \( 1.294\)e\(+00\) & \(3.0\)e\(-07\) & \(0.00\) & \(    15\) & \(    25\) & \(    16\) & first\_order \\
WOODSNE & \(  4000\) & \(  3001\) & \( 5.000\)e\(-01\) & \(0.0\)e\(+00\) & \(0.33\) & \(     2\) & \(     3\) & \(     3\) & first\_order \\\hline
\end{longtable}

\begin{longtable}{lrrrrrrrrl}
\caption{%
    \label{tab:ncl-nls}
    \nclnls results on \(127\) CUTEst nonlinear least-squares problems.
    119 problems were solved successfully.}\\
\hline
name & \(f\)\hspace*{14pt} & \(\|\nabla L\|_2\) & \(\|c\|_2\) & \(t\)\hspace*{6pt} & iter & \#\(c\) & \#\(\nabla c\) & \#\(\nabla^2 L\) & status \\\hline
\endhead
\hline
\multicolumn{10}{r}{{\bfseries Continued on next page}}\\
\hline
\endfoot
\endlastfoot
ARWHDNE & \( 6.971\)e\(+01\) & \(1.5\)e\(-14\) & \(1.4\)e\(-16\) & \(2.65\) & \(    30\) & \(   105\) & \(    32\) & \(    31\) & first\_order \\
BA-L1 & \( 4.400\)e\(-31\) & \(9.1\)e\(-16\) & \(2.0\)e\(-10\) & \(0.01\) & \(     5\) & \(     6\) & \(     7\) & \(     6\) & first\_order \\
BA-L16 & \( 4.324\)e\(+05\) & \(1.4\)e\(+03\) & \(4.5\)e\(-07\) & \(937.15\) & \(    88\) & \(   454\) & \(    90\) & \(    89\) & unknown \\
BA-L1SP & \( 3.204\)e\(-23\) & \(7.4\)e\(-12\) & \(8.3\)e\(-05\) & \(0.01\) & \(     4\) & \(     5\) & \(     6\) & \(     5\) & first\_order \\
BA-L21 & \( 1.975\)e\(+05\) & \(1.5\)e\(-02\) & \(3.0\)e\(-05\) & \(947.20\) & \(   201\) & \(  1093\) & \(   203\) & \(   203\) & max\_time \\
BA-L49 & \( 1.674\)e\(+04\) & \(1.2\)e\(-01\) & \(1.6\)e\(-04\) & \(869.15\) & \(   217\) & \(  1189\) & \(   219\) & \(   219\) & max\_time \\
BA-L52 & \( 3.860\)e\(+06\) & \(2.0\)e\(+02\) & \(8.0\)e\(-01\) & \(1778.45\) & \(    31\) & \(    50\) & \(    33\) & \(    32\) & max\_time \\
BA-L73 & \( 9.609\)e\(+05\) & \(1.1\)e\(+02\) & \(1.9\)e\(-06\) & \(693.06\) & \(   135\) & \(   635\) & \(   137\) & \(   136\) & unknown \\
BARDNE & \( 4.107\)e\(-03\) & \(2.3\)e\(-14\) & \(5.7\)e\(-12\) & \(0.00\) & \(     5\) & \(     6\) & \(     7\) & \(     6\) & first\_order \\
BDQRTICNE & \( 1.000\)e\(+04\) & \(3.4\)e\(-10\) & \(1.7\)e\(-11\) & \(0.47\) & \(    17\) & \(    25\) & \(    19\) & \(    18\) & first\_order \\
BEALENE & \( 4.598\)e\(-20\) & \(2.0\)e\(-10\) & \(1.6\)e\(-09\) & \(0.00\) & \(     8\) & \(     9\) & \(    10\) & \(     9\) & first\_order \\
BIGGS6NE & \( 2.153\)e\(-17\) & \(6.8\)e\(-12\) & \(7.2\)e\(-08\) & \(0.02\) & \(    60\) & \(    74\) & \(    62\) & \(    61\) & first\_order \\
BOX3NE & \( 2.641\)e\(-19\) & \(7.3\)e\(-12\) & \(1.8\)e\(-10\) & \(0.00\) & \(     5\) & \(     6\) & \(     7\) & \(     6\) & first\_order \\
BROWNBSNE & \( 5.933\)e\(-33\) & \(5.6\)e\(-12\) & \(1.1\)e\(-10\) & \(0.00\) & \(    13\) & \(    42\) & \(    15\) & \(    14\) & first\_order \\
BROWNDENE & \( 4.291\)e\(+04\) & \(9.7\)e\(-09\) & \(3.0\)e\(-11\) & \(0.00\) & \(    14\) & \(    22\) & \(    16\) & \(    15\) & first\_order \\
BRYBNDNE & \( 6.662\)e\(-20\) & \(2.6\)e\(-10\) & \(2.5\)e\(-10\) & \(0.20\) & \(     6\) & \(     7\) & \(     8\) & \(     7\) & first\_order \\
CHAINWOONE & \( 6.598\)e\(+03\) & \(1.8\)e\(+00\) & \(2.1\)e\(-01\) & \(14.40\) & \(   500\) & \(   546\) & \(   502\) & \(   501\) & max\_iter \\
CHEBYQADNE & \( 4.358\)e\(-03\) & \(5.1\)e\(-07\) & \(2.9\)e\(-07\) & \(0.60\) & \(    33\) & \(    51\) & \(    35\) & \(    34\) & first\_order \\
CHNRSBNE & \( 2.913\)e\(-21\) & \(1.9\)e\(-10\) & \(7.2\)e\(-11\) & \(0.02\) & \(    34\) & \(    51\) & \(    36\) & \(    35\) & first\_order \\
CHNRSNBMNE & \( 6.544\)e\(-23\) & \(4.0\)e\(-11\) & \(3.9\)e\(-10\) & \(0.02\) & \(    38\) & \(    57\) & \(    40\) & \(    39\) & first\_order \\
COATINGNE & \( 2.527\)e\(-01\) & \(1.5\)e\(-09\) & \(3.9\)e\(-10\) & \(0.02\) & \(    11\) & \(    22\) & \(    13\) & \(    12\) & first\_order \\
CUBENE & \( 2.479\)e\(-33\) & \(9.6\)e\(-17\) & \(6.7\)e\(-15\) & \(0.00\) & \(     2\) & \(     7\) & \(     4\) & \(     3\) & first\_order \\
DECONVBNE & \( 1.285\)e\(-03\) & \(3.1\)e\(-07\) & \(5.3\)e\(-06\) & \(0.03\) & \(    20\) & \(    34\) & \(    23\) & \(    21\) & first\_order \\
DECONVNE & \( 2.694\)e\(-15\) & \(3.8\)e\(-08\) & \(6.4\)e\(-10\) & \(0.00\) & \(     2\) & \(     3\) & \(     4\) & \(     3\) & first\_order \\
DENSCHNBNE & \( 1.225\)e\(-17\) & \(2.6\)e\(-09\) & \(5.8\)e\(-09\) & \(0.00\) & \(     6\) & \(     8\) & \(     8\) & \(     7\) & first\_order \\
DENSCHNCNE & \( 9.992\)e\(-38\) & \(4.3\)e\(-19\) & \(6.7\)e\(-06\) & \(0.00\) & \(     6\) & \(     7\) & \(     8\) & \(     7\) & first\_order \\
DENSCHNDNE & \( 4.012\)e\(-28\) & \(4.5\)e\(-17\) & \(3.6\)e\(-04\) & \(0.00\) & \(    15\) & \(    16\) & \(    17\) & \(    16\) & first\_order \\
DENSCHNENE & \( 4.540\)e\(-21\) & \(9.5\)e\(-11\) & \(1.7\)e\(-06\) & \(0.00\) & \(    15\) & \(    20\) & \(    17\) & \(    16\) & first\_order \\
DENSCHNFNE & \( 2.351\)e\(-38\) & \(2.2\)e\(-19\) & \(1.9\)e\(-06\) & \(0.00\) & \(     4\) & \(     5\) & \(     6\) & \(     5\) & first\_order \\
DEVGLA1NE & \( 1.063\)e\(-13\) & \(2.9\)e\(-08\) & \(4.0\)e\(-11\) & \(0.01\) & \(    16\) & \(    45\) & \(    18\) & \(    17\) & first\_order \\
DEVGLA2NE & \( 1.672\)e\(-14\) & \(6.3\)e\(-08\) & \(3.0\)e\(-07\) & \(0.00\) & \(     9\) & \(    12\) & \(    11\) & \(    10\) & first\_order \\
EGGCRATENE & \( 4.744\)e\(+00\) & \(3.1\)e\(-08\) & \(2.4\)e\(-09\) & \(0.00\) & \(     4\) & \(     5\) & \(     6\) & \(     5\) & first\_order \\
ELATVIDUNE & \( 2.738\)e\(+01\) & \(1.2\)e\(-09\) & \(5.0\)e\(-10\) & \(0.00\) & \(     8\) & \(     9\) & \(    10\) & \(     9\) & first\_order \\
ENGVAL2NE & \( 6.374\)e\(-15\) & \(1.3\)e\(-08\) & \(1.2\)e\(-07\) & \(0.00\) & \(    12\) & \(    20\) & \(    14\) & \(    13\) & first\_order \\
ERRINROSNE & \( 2.020\)e\(+01\) & \(1.4\)e\(-09\) & \(4.1\)e\(-10\) & \(0.01\) & \(    17\) & \(    24\) & \(    19\) & \(    18\) & first\_order \\
ERRINRSMNE & \( 1.926\)e\(+01\) & \(1.7\)e\(-08\) & \(1.2\)e\(-09\) & \(0.02\) & \(    25\) & \(    38\) & \(    27\) & \(    26\) & first\_order \\
EXP2NE & \( 1.617\)e\(-19\) & \(7.1\)e\(-13\) & \(1.9\)e\(-13\) & \(0.00\) & \(     5\) & \(     6\) & \(     7\) & \(     6\) & first\_order \\
EXPFITNE & \( 1.203\)e\(-01\) & \(2.4\)e\(-10\) & \(4.9\)e\(-12\) & \(0.01\) & \(    21\) & \(   149\) & \(    23\) & \(    22\) & first\_order \\
EXTROSNBNE & \(-2.002\)e\(+00\) & \(1.0\)e\(-06\) & \(1.3\)e\(-09\) & \(1.16\) & \(   250\) & \(  1860\) & \(   252\) & \(   251\) & first\_order \\
FBRAIN2NE & \( 1.842\)e\(-01\) & \(6.3\)e\(-10\) & \(3.1\)e\(-11\) & \(0.11\) & \(     5\) & \(     6\) & \(     7\) & \(     6\) & first\_order \\
FBRAINNE & \( 2.083\)e\(-01\) & \(8.4\)e\(-08\) & \(1.7\)e\(-09\) & \(0.05\) & \(     4\) & \(     5\) & \(     6\) & \(     5\) & first\_order \\
FREURONE & \( 2.449\)e\(+01\) & \(1.5\)e\(-10\) & \(5.5\)e\(-12\) & \(0.00\) & \(     7\) & \(    15\) & \(     9\) & \(     8\) & first\_order \\
GENROSEBNE & \( 7.965\)e\(+02\) & \(1.3\)e\(-06\) & \(1.5\)e\(-06\) & \(0.04\) & \(     6\) & \(     8\) & \(     9\) & \(     7\) & first\_order \\
GENROSENE & \( 5.000\)e\(-01\) & \(8.9\)e\(-14\) & \(1.1\)e\(-14\) & \(2.07\) & \(   436\) & \(   639\) & \(   438\) & \(   437\) & first\_order \\
GULFNE & \( 1.755\)e\(-20\) & \(3.4\)e\(-11\) & \(1.6\)e\(-10\) & \(0.02\) & \(    20\) & \(    28\) & \(    22\) & \(    21\) & first\_order \\
HATFLDANE & \( 8.445\)e\(-20\) & \(5.4\)e\(-11\) & \(6.2\)e\(-09\) & \(0.00\) & \(     9\) & \(    10\) & \(    11\) & \(    10\) & first\_order \\
HATFLDBNE & \( 2.786\)e\(-03\) & \(4.8\)e\(-12\) & \(1.1\)e\(-11\) & \(0.00\) & \(     7\) & \(     8\) & \(     9\) & \(     8\) & first\_order \\
HATFLDCNE & \( 2.674\)e\(-17\) & \(4.6\)e\(-09\) & \(7.8\)e\(-09\) & \(0.00\) & \(     3\) & \(     4\) & \(     5\) & \(     4\) & first\_order \\
HATFLDDNE & \( 1.273\)e\(-07\) & \(4.1\)e\(-09\) & \(9.7\)e\(-07\) & \(0.00\) & \(     5\) & \(     8\) & \(     7\) & \(     6\) & first\_order \\
HATFLDENE & \( 1.364\)e\(-06\) & \(5.8\)e\(-11\) & \(1.7\)e\(-09\) & \(0.00\) & \(     5\) & \(     6\) & \(     7\) & \(     6\) & first\_order \\
HATFLDFLNE & \( 3.008\)e\(-05\) & \(9.2\)e\(-12\) & \(7.3\)e\(-09\) & \(0.02\) & \(   103\) & \(   443\) & \(   105\) & \(   104\) & first\_order \\
HELIXNE & \( 5.661\)e\(-43\) & \(1.7\)e\(-21\) & \(1.7\)e\(-10\) & \(0.00\) & \(     9\) & \(    12\) & \(    11\) & \(    10\) & first\_order \\
HIMMELBFNE & \( 2.168\)e\(+06\) & \(2.1\)e\(+02\) & \(2.9\)e\(-04\) & \(0.09\) & \(   500\) & \(   720\) & \(   502\) & \(   501\) & max\_iter \\
HS1NE & \( 3.359\)e\(-17\) & \(4.1\)e\(-09\) & \(5.8\)e\(-16\) & \(0.00\) & \(     7\) & \(    11\) & \(     9\) & \(     8\) & first\_order \\
HS25NE & \( 1.642\)e\(+01\) & \(8.1\)e\(-08\) & \(3.6\)e\(-11\) & \(0.00\) & \(     3\) & \(     5\) & \(     6\) & \(     4\) & first\_order \\
HS2NE & \( 2.471\)e\(+00\) & \(8.4\)e\(-12\) & \(1.3\)e\(-12\) & \(0.00\) & \(     5\) & \(     7\) & \(     8\) & \(     6\) & first\_order \\
INTEQNE & \( 6.816\)e\(-37\) & \(8.7\)e\(-19\) & \(1.6\)e\(-07\) & \(0.00\) & \(     2\) & \(     3\) & \(     4\) & \(     3\) & first\_order \\
JENSMPNE & \( 6.218\)e\(+01\) & \(1.6\)e\(-08\) & \(3.1\)e\(-09\) & \(0.00\) & \(     8\) & \(    14\) & \(    10\) & \(     9\) & first\_order \\
JUDGENE & \( 8.041\)e\(+00\) & \(4.8\)e\(-07\) & \(2.5\)e\(-07\) & \(0.00\) & \(     5\) & \(     6\) & \(     7\) & \(     6\) & first\_order \\
KOEBHELBNE & \( 3.876\)e\(+01\) & \(1.4\)e\(-08\) & \(4.5\)e\(-09\) & \(0.13\) & \(   119\) & \(   164\) & \(   121\) & \(   121\) & first\_order \\
KOWOSBNE & \( 1.539\)e\(-04\) & \(2.9\)e\(-10\) & \(2.6\)e\(-09\) & \(0.00\) & \(    11\) & \(    14\) & \(    13\) & \(    12\) & first\_order \\
LIARWHDNE & \( 1.113\)e\(-21\) & \(4.7\)e\(-11\) & \(4.0\)e\(-11\) & \(0.17\) & \(     6\) & \(     7\) & \(     8\) & \(     7\) & first\_order \\
LINVERSENE & \( 3.405\)e\(+02\) & \(3.3\)e\(-07\) & \(1.4\)e\(-07\) & \(0.31\) & \(    19\) & \(    22\) & \(    22\) & \(    20\) & first\_order \\
MANCINONE & \( 2.507\)e\(-32\) & \(1.1\)e\(-16\) & \(2.4\)e\(-09\) & \(0.12\) & \(     4\) & \(     5\) & \(     6\) & \(     5\) & first\_order \\
MANNE & \(-9.698\)e\(-01\) & \(9.6\)e\(-07\) & \(0.0\)e\(+00\) & \(10.27\) & \(   344\) & \(   346\) & \(   347\) & \(   345\) & first\_order \\
MARINE & \( 2.009\)e\(+06\) & \(1.7\)e\(-08\) & \(1.1\)e\(-08\) & \(3.76\) & \(    40\) & \(    43\) & \(    43\) & \(    41\) & first\_order \\
MEYER3NE & \( 4.397\)e\(+01\) & \(2.0\)e\(-06\) & \(1.1\)e\(-08\) & \(0.00\) & \(     8\) & \(    14\) & \(    10\) & \(     9\) & first\_order \\
MODBEALENE & \( 2.521\)e\(-18\) & \(1.5\)e\(-09\) & \(2.8\)e\(-09\) & \(1.80\) & \(    11\) & \(    12\) & \(    13\) & \(    12\) & first\_order \\
MOREBVNE & \( 2.674\)e\(-37\) & \(9.3\)e\(-19\) & \(7.6\)e\(-08\) & \(0.00\) & \(     2\) & \(     3\) & \(     4\) & \(     3\) & first\_order \\
MUONSINE & \( 2.194\)e\(+04\) & \(2.4\)e\(-11\) & \(6.2\)e\(-14\) & \(0.02\) & \(    10\) & \(    15\) & \(    12\) & \(    11\) & first\_order \\
NGONE & \(-8.551\)e\(+00\) & \(4.0\)e\(-09\) & \(3.8\)e\(-11\) & \(2.46\) & \(   106\) & \(   108\) & \(   109\) & \(   107\) & first\_order \\
NONDIANE & \( 4.949\)e\(-01\) & \(1.3\)e\(-16\) & \(1.5\)e\(-07\) & \(0.16\) & \(     6\) & \(     7\) & \(     8\) & \(     7\) & first\_order \\
NONMSQRTNE & \( 3.543\)e\(+02\) & \(1.1\)e\(-07\) & \(2.5\)e\(-06\) & \(4.01\) & \(    19\) & \(    30\) & \(    21\) & \(    20\) & first\_order \\
NONSCOMPNE & \( 1.805\)e\(-07\) & \(1.8\)e\(-10\) & \(3.0\)e\(-06\) & \(0.52\) & \(    19\) & \(    28\) & \(    21\) & \(    20\) & first\_order \\
OSCIGRNE & \( 1.449\)e\(-35\) & \(4.0\)e\(-18\) & \(1.4\)e\(-07\) & \(3.78\) & \(     6\) & \(     7\) & \(     8\) & \(     7\) & first\_order \\
OSCIPANE & \( 5.000\)e\(-01\) & \(9.8\)e\(-06\) & \(5.5\)e\(-04\) & \(0.16\) & \(   500\) & \(  2698\) & \(   502\) & \(   501\) & max\_iter \\
PALMER1ANE & \( 2.988\)e\(+01\) & \(2.4\)e\(-06\) & \(6.8\)e\(-08\) & \(0.01\) & \(    12\) & \(    16\) & \(    14\) & \(    13\) & first\_order \\
PALMER1BNE & \( 1.724\)e\(+00\) & \(1.9\)e\(-08\) & \(1.9\)e\(-10\) & \(0.01\) & \(    11\) & \(    16\) & \(    13\) & \(    12\) & first\_order \\
PALMER1ENE & \( 4.176\)e\(-04\) & \(7.7\)e\(-09\) & \(1.4\)e\(-06\) & \(0.00\) & \(     7\) & \(     8\) & \(     9\) & \(     8\) & first\_order \\
PALMER1NE & \( 5.877\)e\(+03\) & \(2.7\)e\(-06\) & \(2.4\)e\(-09\) & \(0.02\) & \(    37\) & \(    46\) & \(    39\) & \(    38\) & first\_order \\
PALMER2ANE & \( 8.555\)e\(-03\) & \(1.9\)e\(-10\) & \(4.1\)e\(-08\) & \(0.02\) & \(    44\) & \(    69\) & \(    46\) & \(    45\) & first\_order \\
PALMER2BNE & \( 3.116\)e\(-01\) & \(1.9\)e\(-07\) & \(3.1\)e\(-09\) & \(0.00\) & \(     8\) & \(    11\) & \(    10\) & \(     9\) & first\_order \\
PALMER2ENE & \( 5.815\)e\(-02\) & \(6.9\)e\(-08\) & \(3.2\)e\(-10\) & \(0.01\) & \(    20\) & \(    23\) & \(    22\) & \(    21\) & first\_order \\
PALMER2NE & \( 1.826\)e\(+03\) & \(3.4\)e\(-06\) & \(5.3\)e\(-08\) & \(0.02\) & \(    49\) & \(    64\) & \(    51\) & \(    50\) & first\_order \\
PALMER3ANE & \( 1.022\)e\(-02\) & \(2.0\)e\(-08\) & \(9.0\)e\(-07\) & \(0.01\) & \(    13\) & \(    17\) & \(    15\) & \(    14\) & first\_order \\
PALMER3BNE & \( 2.114\)e\(+00\) & \(6.1\)e\(-12\) & \(1.1\)e\(-14\) & \(0.01\) & \(    27\) & \(    38\) & \(    29\) & \(    28\) & first\_order \\
PALMER3ENE & \( 2.537\)e\(-05\) & \(4.1\)e\(-11\) & \(1.2\)e\(-09\) & \(0.01\) & \(    14\) & \(    19\) & \(    16\) & \(    15\) & first\_order \\
PALMER3NE & \( 1.208\)e\(+03\) & \(1.1\)e\(-05\) & \(2.4\)e\(-08\) & \(0.00\) & \(     6\) & \(    10\) & \(     8\) & \(     7\) & first\_order \\
PALMER4ANE & \( 2.030\)e\(-02\) & \(2.8\)e\(-08\) & \(1.6\)e\(-07\) & \(0.01\) & \(    11\) & \(    20\) & \(    13\) & \(    12\) & first\_order \\
PALMER4BNE & \( 3.418\)e\(+00\) & \(1.3\)e\(-12\) & \(2.3\)e\(-14\) & \(0.01\) & \(    32\) & \(    41\) & \(    34\) & \(    33\) & first\_order \\
PALMER4ENE & \( 7.400\)e\(-05\) & \(4.7\)e\(-09\) & \(2.0\)e\(-07\) & \(0.01\) & \(    11\) & \(    13\) & \(    13\) & \(    12\) & first\_order \\
PALMER4NE & \( 1.212\)e\(+03\) & \(1.5\)e\(-05\) & \(4.2\)e\(-08\) & \(0.00\) & \(     6\) & \(    10\) & \(     8\) & \(     7\) & first\_order \\
PALMER5ANE & \( 1.452\)e\(-02\) & \(9.5\)e\(-07\) & \(4.3\)e\(-07\) & \(0.10\) & \(   182\) & \(   874\) & \(   184\) & \(   183\) & first\_order \\
PALMER5BNE & \( 4.876\)e\(-03\) & \(3.7\)e\(-09\) & \(1.1\)e\(-07\) & \(0.01\) & \(    44\) & \(    50\) & \(    46\) & \(    45\) & first\_order \\
PALMER5ENE & \( 1.036\)e\(-02\) & \(1.3\)e\(-09\) & \(3.3\)e\(-07\) & \(0.02\) & \(    60\) & \(   200\) & \(    62\) & \(    61\) & first\_order \\
PALMER6ANE & \( 2.797\)e\(-02\) & \(2.1\)e\(-08\) & \(7.5\)e\(-06\) & \(0.01\) & \(    22\) & \(    31\) & \(    24\) & \(    23\) & first\_order \\
PALMER6ENE & \( 6.462\)e\(-02\) & \(5.5\)e\(-07\) & \(2.3\)e\(-06\) & \(0.00\) & \(     7\) & \(    10\) & \(     9\) & \(     8\) & first\_order \\
PALMER7ANE & \( 5.168\)e\(+00\) & \(1.6\)e\(-07\) & \(2.6\)e\(-05\) & \(0.16\) & \(   355\) & \(  1551\) & \(   357\) & \(   358\) & first\_order \\
PALMER7ENE & \( 3.346\)e\(+00\) & \(1.5\)e\(-06\) & \(4.0\)e\(-10\) & \(0.03\) & \(    54\) & \(   274\) & \(    56\) & \(    55\) & first\_order \\
PALMER8ANE & \( 3.700\)e\(-02\) & \(9.3\)e\(-08\) & \(3.6\)e\(-07\) & \(0.01\) & \(    19\) & \(    31\) & \(    21\) & \(    20\) & first\_order \\
PALMER8ENE & \( 3.170\)e\(-01\) & \(6.3\)e\(-08\) & \(3.3\)e\(-09\) & \(0.01\) & \(    25\) & \(    35\) & \(    27\) & \(    26\) & first\_order \\
PENLT1NE & \( 3.544\)e\(-10\) & \(5.5\)e\(-14\) & \(3.8\)e\(-05\) & \(0.00\) & \(     8\) & \(     9\) & \(    10\) & \(     9\) & first\_order \\
PENLT2NE & \( 4.688\)e\(-11\) & \(1.4\)e\(-18\) & \(5.6\)e\(-08\) & \(0.00\) & \(     5\) & \(    10\) & \(     7\) & \(     6\) & first\_order \\
PINENE & \( 4.194\)e\(-05\) & \(4.9\)e\(-10\) & \(1.1\)e\(-06\) & \(0.72\) & \(    13\) & \(    18\) & \(    16\) & \(    14\) & first\_order \\
POWERSUMNE & \( 1.799\)e\(-22\) & \(1.2\)e\(-11\) & \(6.0\)e\(-10\) & \(0.01\) & \(    45\) & \(    48\) & \(    47\) & \(    46\) & first\_order \\
PRICE3NE & \( 1.505\)e\(-35\) & \(5.2\)e\(-18\) & \(1.1\)e\(-05\) & \(0.00\) & \(     6\) & \(     7\) & \(     8\) & \(     7\) & first\_order \\
PRICE4NE & \( 3.508\)e\(-32\) & \(3.4\)e\(-19\) & \(1.1\)e\(-04\) & \(0.00\) & \(    13\) & \(    14\) & \(    15\) & \(    14\) & first\_order \\
QINGNE & \( 1.318\)e\(-33\) & \(2.8\)e\(-17\) & \(9.2\)e\(-05\) & \(0.00\) & \(     4\) & \(     7\) & \(     6\) & \(     5\) & first\_order \\
RSNBRNE & \( 5.556\)e\(-32\) & \(3.3\)e\(-16\) & \(3.1\)e\(-14\) & \(0.00\) & \(     1\) & \(     3\) & \(     3\) & \(     2\) & first\_order \\
S308NE & \( 3.866\)e\(-01\) & \(1.7\)e\(-09\) & \(9.0\)e\(-10\) & \(0.00\) & \(    10\) & \(    16\) & \(    12\) & \(    11\) & first\_order \\
SBRYBNDNE & \( 6.690\)e\(-20\) & \(2.6\)e\(-10\) & \(2.5\)e\(-10\) & \(0.18\) & \(     6\) & \(     7\) & \(     8\) & \(     7\) & first\_order \\
SINVALNE & \( 1.578\)e\(-30\) & \(1.8\)e\(-15\) & \(1.8\)e\(-15\) & \(0.00\) & \(     1\) & \(     3\) & \(     3\) & \(     2\) & first\_order \\
SPECANNE & \( 3.291\)e\(-13\) & \(3.5\)e\(-08\) & \(1.2\)e\(-12\) & \(0.49\) & \(    10\) & \(    14\) & \(    12\) & \(    11\) & first\_order \\
SROSENBRNE & \( 5.892\)e\(-33\) & \(1.1\)e\(-18\) & \(2.2\)e\(-18\) & \(0.05\) & \(     2\) & \(     3\) & \(     4\) & \(     3\) & first\_order \\
SSBRYBNDNE & \( 6.663\)e\(-20\) & \(2.6\)e\(-10\) & \(2.5\)e\(-10\) & \(0.21\) & \(     6\) & \(     7\) & \(     8\) & \(     7\) & first\_order \\
STREGNE & \( 2.204\)e\(-27\) & \(6.6\)e\(-22\) & \(4.4\)e\(-15\) & \(0.00\) & \(     2\) & \(     3\) & \(     4\) & \(     3\) & first\_order \\
STRTCHDVNE & \( 2.260\)e\(-15\) & \(1.8\)e\(-10\) & \(6.1\)e\(-07\) & \(0.00\) & \(    10\) & \(    11\) & \(    12\) & \(    11\) & first\_order \\
TQUARTICNE & \( 0.000\)e\(+00\) & \(0.0\)e\(+00\) & \(2.2\)e\(-16\) & \(0.05\) & \(     1\) & \(     2\) & \(     3\) & \(     2\) & first\_order \\
TRIGON1NE & \( 3.427\)e\(-39\) & \(5.7\)e\(-20\) & \(1.8\)e\(-08\) & \(0.00\) & \(     4\) & \(     5\) & \(     6\) & \(     5\) & first\_order \\
TRIGON2NE & \( 1.620\)e\(+00\) & \(2.7\)e\(-09\) & \(1.3\)e\(-09\) & \(0.00\) & \(    11\) & \(    12\) & \(    13\) & \(    12\) & first\_order \\
VARDIMNE & \( 4.873\)e\(-20\) & \(2.9\)e\(-09\) & \(1.2\)e\(-10\) & \(0.00\) & \(    14\) & \(    19\) & \(    16\) & \(    15\) & first\_order \\
VIBRBEAMNE & \( 7.822\)e\(-02\) & \(2.3\)e\(-11\) & \(1.7\)e\(-14\) & \(0.01\) & \(     7\) & \(     8\) & \(     9\) & \(     8\) & first\_order \\
WATSONNE & \( 1.430\)e\(-15\) & \(5.8\)e\(-16\) & \(5.6\)e\(-14\) & \(0.00\) & \(     4\) & \(     5\) & \(     6\) & \(     5\) & first\_order \\
WAYSEA1NE & \( 0.000\)e\(+00\) & \(0.0\)e\(+00\) & \(2.7\)e\(-08\) & \(0.00\) & \(     7\) & \(     8\) & \(     9\) & \(     8\) & first\_order \\
WAYSEA2NE & \( 8.386\)e\(-17\) & \(2.3\)e\(-09\) & \(1.0\)e\(-07\) & \(0.00\) & \(    10\) & \(    20\) & \(    12\) & \(    11\) & first\_order \\
WEEDSNE & \( 1.294\)e\(+00\) & \(8.6\)e\(-08\) & \(7.4\)e\(-07\) & \(0.01\) & \(    17\) & \(    24\) & \(    19\) & \(    18\) & first\_order \\
WOODSNE & \(-1.906\)e\(+04\) & \(9.1\)e\(-13\) & \(8.0\)e\(-10\) & \(0.04\) & \(     3\) & \(     4\) & \(     5\) & \(     4\) & first\_order \\\hline
\end{longtable}

\end{document}